\documentclass[12pt]{article}
\usepackage{amssymb,amsfonts,color}
\usepackage{amscd}
\usepackage{verbatim}

\newtheorem{theorem}{Theorem}[section]
\newtheorem{lemma}[theorem]{Lemma}
\newtheorem{proposition}[theorem]{Proposition}

\newtheorem{definition}[theorem]{Definition}
\newtheorem{corollary}[theorem]{Corollary}
\newtheorem{remark}[theorem]{Remark}

\newenvironment{proof}{\bf Proof. \rm}{$\Box$}

\newcommand{\be}{\begin{equation}}
\newcommand{\ee}{\end{equation}}

\newcommand{\norm}[1]{\Vert #1 \Vert}
\newcommand{\cA}{\mathcal{A}}

\newcommand{\cF}{\mathcal{F}}
\newcommand{\cL}{\mathcal{L}}
\newcommand{\cO}{\mathcal{O}}

\newcommand{\bB}{{\mathbb{B}}}
\newcommand{\bC}{{\mathbb{C}}}

\newcommand{\bZ}{{\mathbb{Z}}}
\newcommand{\Ftheta}{\mathbb{F}^+_\theta}

\begin{document}

\title{Operator Algebras Associated with Unitary Commutation
Relations\renewcommand{\thefootnote}{}
}
\author{}
\author{
~\\ Stephen C. Power\thanks{SCP is supported by EPSRC grant
EP/E002625/1 }\\[1ex] {\small\itshape Lancaster University }\\
{\small\itshape Department of Mathematics and Statistics }\\
{\small\itshape Lancaster, United Kingdom LA1 4YF}\\
{\small{\itshape E-mail: \tt{s.power@lancaster.ac.uk}}}\\ ~\\
Baruch Solel\thanks{BS is supported by the Fund for the Promotion
of Research at the Technion and by EPSRC grant
EP/E002625/1}\\[1ex] {\small\itshape Technion}\\ {\small\itshape
Department of Mathematics}\\ {\small\itshape Haifa 32000,
Israel}\\ {\small{\itshape E-mail:
\tt{mabaruch@techunix.technion.ac.il}}}\\
\\
}

\renewcommand{\baselinestretch}{1}
\maketitle

%\maketitle

\newpage

\begin{abstract}
We define nonselfadjoint operator algebras with generators\\
$L_{e_1},\dots, L_{e_n}, L_{f_1},\dots,L_{f_m}$ subject to the
unitary commutation relations of the form
\[ L_{e_i}L_{f_j} = \sum_{k,l} u_{i,j,k,l} L_{f_l}L_{e_k}
\]
where $u= (u_{i,j,k,l})$ is an $nm \times nm$ unitary matrix.
These algebras, which generalise the analytic Toeplitz algebras of
rank 2 graphs with a single vertex, are classified up to isometric
isomorphism in terms of the matrix $u$.
\end{abstract}

\begin{section}{Introduction}
The unilateral shift on complex separable Hilbert space generates
two fundamental operator algebras, namely the norm closed (unital)
algebra
%, which is naturally isomorphic to the disc algebra, and
and the weak operator topology  closed algebra. The former is
naturally isomorphic to the disc algebra of holomorphic functions
on the unit disc, continuous to the boundary, while the latter is
isomorphic to  $H^\infty $. The freely noncommuting multivariable
generalisations of these algebras arise from the freely
noncommuting shifts $L_{e_1},\dots , L_{e_n}$ given by the left
creation operators on the Fock space $\cF_n = \sum_{k=0}^\infty
\oplus (\bC^n)^{\otimes k}.$ Here the generated operator algebras,
denoted $\cA_n$ and $\cL_n$ for the norm and weak topologies, are
known as the noncommutative disc algebra and the freesemigroup
algebra. They  have been studied extensively  with respect to
operator algebra structure, representation theory and the
multivariable operator theory of row contractions. See for example
 \cite{DP98}, \cite{Pop}.

Higher rank generalisations of these algebras arise when one
considers  several families of freely noncommuting generators
between which there are commutation relations. In the present
paper we consider a very general form of such relations, namely
\[ L_{e_i}L_{f_j} = \sum_{k,l} u_{i,j,k,l} L_{f_l}L_{e_k}
\]
where $L_{e_1},\dots , L_{e_n}$ and $L_{f_1},\dots , L_{f_m}$ are
freely noncommuting and $u= (u_{i,j,k,l})$ is an $nm \times nm$
unitary matrix. The associated operator algebras are denoted
$\cA_u$ and $\cL_u$ and we classify them up to various forms of
isomorphism in terms of the unitary matrices $u$. Such unitary
relations arose originally in the context of the general dilation
theorem proven in Solel (\cite{S}, \cite{S2}) for two row
contractions $[T_1 \cdots T_n]$ and $ [S_1 \cdots S_m]$ satisfying
the unitary commutation relations.

 For $n=m=1$, we have $u = [\alpha]$ with $|\alpha| = 1$ and $\cA_u$ is the
subalgebra of the rotation C*-algebra for the relations $uv
=\alpha vu$. When $u$ is a permutation unitary matrix arising from
a permutation $\theta$ in $S_{nm}$ then the relations are those
associated with a single vertex rank 2 graph in the sense of
Kumjian and Pask, and the algebras in this case have been
considered in Kribs and Power \cite{KP} and Power \cite{Po}. In
particular, in \cite{Po} it was shown that there are 9 operator
algebras $\cA_\theta$ arising from the 24 permutations in case
$n=m=2$. In contrast, we see below in Section 6 that for general
$2$ by $2$ unitaries $u$  there are uncountably many isomorphism
classes of the unitary relation algebras $\cA_u$ expressed in
terms of a nine fold real parametrisation of isomorphism types.

The algebras $\cA_\theta$ are easily defined; they are determined
by the left regular representation of the semigroup $\Ftheta$
whose generators are $e_1,\dots ,e_n, f_1,\dots ,f_m$ subject to
the relations $e_if_j = f_le_k$ where $\theta(i,j)=(k,l)$. On the
other hand the unitary relation algebras $\cA_u$ are generated by
creation operators on a $\bZ_+^2$-graded Fock space $ \sum_{k,l}
\oplus (\bC^n)^{\otimes k}\otimes (\bC^m)^{\otimes l}$ with
relations arising from the identification $u : \bC^n \otimes \bC^m
\to \bC^m \otimes \bC^n$. In particular, $\cA_u$ is a
representation of the non-selfadjoint tensor algebra of a
 rank $2$ correspondence (or a product system over
$\mathbb{N}^2$) in the sense of \cite{S2}. See also \cite{F}

In the main results, summarised partly in Theorem
\ref{bigradedisom}, we see that if $ \mathcal{A}_u $ and $
\mathcal{A}_v$ are isomorphic then  the two families of generators
have matching cardinalities. Furthermore, if $n\neq m$ then the
algebras are isomorphic if and only if the unitaries $u, v$ in
$M_{nm}(\bC)$ are unitary equivalent by a unitary $A\otimes B$ in
$M_n(\bC) \otimes M_m(\bC)$. As in \cite{Po} we term this
\emph{product unitary equivalence} (with respect to the fixed
tensor product decomposition). The case $n=m$ admits an extra
possibility, in view of the possibility of generator exchanging
isomorphisms, namely that
 $u, \tilde{v}$
are product unitary equivalent, where
$\tilde{v}_{i,j,k,l}=\bar{v}_{l,k,j,i}$.

The theorem is proven as follows. After some preliminaries we
identify, in Section 3,  the character space $M(\cA_u)$ and the
set of w*-continuous characters on $\cL_u$. These are  subsets of
the closed unit ball product $\overline{\bB}_n \times
\overline{\bB}_m$ which are associated with a variety $V_u$ in
$\bC^n \times \bC^m$ determined by $u$. We then define the
\emph{core}  $\Omega_u^0$, a closed subset of the realised
character space $ \Omega_u = M(\cA_u)$, and we identify this
intrinsically (algebraically) in terms of representations of
$\cA_u$ into $T_2$, the algebra of upper triangular matrices in
$M_2(\bC)$. The importance of the core is that we are able to show
that the interior is a minimal automorphism invariant subset on
which automorphisms act transitively. This allows us to infer the
existence of graded isomorphisms from general isomorphisms. To
construct automorphisms we first review, in Section 4,
Voiculescu's construction  of a unitary action of the Lie group
$U(1,n)$ on the Cuntz algebra $\cO_n$ and the operator algebras
$\cA_n$ and $ \cL_n$. This provides, in particular, unitary
automorphisms $\Theta_\alpha$, for $\alpha \in \bB_n$, which act
transitively on the interior ball, $ \bB_n$, of the character
space of $\cA_n$. For these explicit unitary automorphisms of the
$e_i$-generated copy of $\cA_n$ in $\cA_u$, we establish unitary
commutation relations for the tuples $\Theta_\alpha(L_{e_1}),
\dots , \Theta_\alpha(L_{e_n})$ and $L_{f_1},\dots ,L_{f_m}$, when
$(\alpha,0)$ is a point in the core. This enables us to define
natural unitary automorphisms of $\cA_u$ itself, and in Theorem
\ref{corechar} the relative interior of the core is identified as
an automorphism invariant set in the Gelfand space $\Omega_u$. In
Section 5 we determine the graded and bigraded isomorphisms in
terms of product unitary equivalence. To do this we observe that
such isomorphisms induce an origin preserving biholomorphic map
between the cores $\Omega_u^0$ and $\Omega_v^0$ and that these
maps, by a generalised Schwarz's Lemma, are implemented by a
product unitary. We then prove the main classification theorem.

In Section 6 we analyse in detail the case $n=m=2$ and consider
the special case of permutation unitaries.

Finally, in Section 7 we show that the algebra $\mathcal{A}_u$ is
contained in a tensor algebra $\mathcal{T}_+(X)$, associated with
a correspondence $X$ as in \cite{MS98}. Moreover, at least when
$n\neq m$, every automorphism of $\mathcal{A}_u$ extends to an
automorphism of $\mathcal{T}_+(X)$. The advantage of the tensor
algebra is that its representation theory is known (\cite{MS98})
while this is not the case yet for the algebra $\mathcal{A}_u$.
\end{section}

\begin{section}{Preliminaries}
Fix two finite dimensional Hilbert spaces $E=\mathbb{C}^n$ and
$F=\mathbb{C}^m$ and a unitary $mn \times mn$ matrix $u$. The rows
and columns of $u$ are indexed by $\{1,\ldots,n\}\times
\{1,\ldots, m\}$ ($u=(u_{(i,j),(k,l)})$) and when we write $u$ as
an $mn\times mn$ matrix we assume that $\{1,\ldots,n\}\times
\{1,\ldots, m\}$ is ordered lexicographically (so that, for
example, the second row is the row indexed by $(1,2)$). We also
fix orthonormal bases $\{e_i\}$ and $\{f_j\}$ for $E$ and $F$
respectively and the matrix $u$ is used to identify $E\otimes F$
with $F\otimes E$ through the equation
\begin{equation}\label{relations}
e_i \otimes f_j = \sum_{k,l} u_{(i,j),(k,l)} f_l \otimes e_k .
\end{equation}
Equivalently, we write \begin{equation}\label{relations2}
f_l\otimes e_k=\sum_{i,j} \bar{u}_{(i,j),(k,l)}e_i\otimes
f_j.\end{equation}
 For every $k,l \in
\mathbb{N}$, we write $X(k,l)$ for $E^{\otimes k}\otimes
F^{\otimes l}$. Using succesive applications of (\ref{relations}),
we can identify $X(k,l)$ with $E^{\otimes k_1}\otimes F^{\otimes
l_1}\otimes E^{\otimes k_2} \otimes \cdots \otimes F^{\otimes
l_r}$ whenever $k=\sum k_i$ and $l=\sum l_j$.

Let $\mathcal{F}(n,m,u)$ be the Fock space given by the Hilbert
space direct sum
$$ \sum_{k,l}X(k,l)=\sum_{k,l} E^{\otimes k}\otimes F^{\otimes l},
$$ and, for $e\in E$ and $f\in F$, write $L_e$ and $L_f$ for the
``shift" operators $$ L_e e_{i_1}\otimes e_{i_2} \otimes \cdots
\otimes e_{i_k}\otimes f_{j_1}\otimes f_{j_2} \otimes \cdots
\otimes f_{j_l} =e \otimes e_{i_1}\otimes e_{i_2} \otimes \cdots
\otimes e_{i_k}\otimes f_{j_1}\otimes f_{j_2} \otimes \cdots
\otimes f_{j_l} $$ and $$ L_f e_{i_1}\otimes e_{i_2} \otimes
\cdots \otimes e_{i_k}\otimes f_{j_1}\otimes f_{j_2} \otimes
\cdots \otimes f_{j_l} =f \otimes e_{i_1}\otimes e_{i_2} \otimes
\cdots \otimes e_{i_k}\otimes f_{j_1}\otimes f_{j_2} \otimes
\cdots \otimes f_{j_l}$$ where, in the last equation, we use
(\ref{relations}) to identify the resulting vector as a vector of
$E^{\otimes k}\otimes F^{\otimes (l+1)}$.

The unital semigroup generated by $\{I,L_e,L_f \;:\; e\in E,\;f\in
F\}$ is denoted $\mathbb{F}^+_u$ and the algebra it generates
denoted $\mathbb{C}[\mathbb{F}_u^+]$. The norm closure of
$\mathbb{C}[\mathbb{F}_u^+]$ will be written $\mathcal{A}_u$ and
its closure in the weak* operator topology will be written
$\mathcal{L}_u$. In particular, the algebras
$\mathcal{L}_{\theta}$ and $\mathcal{A}_{\theta}$ studied in
\cite{Po} are the algebras $\mathcal{L}_u$ and $\mathcal{A}_u$ for
$u$ which is a permutation matrix.

The results of Section 2 in \cite{KP} hold here too with minor
changes. Every $A\in \mathcal{L}_u$ is the limit (in the strong
operator topology) of its Cesaro sums $$\Sigma_p(A)=\sum_{k\leq
p}(1-\frac{k}{p})\Phi_k(A)$$ where $\Phi_k(A)$ lies in
$\mathcal{L}_u$ and is ``supported" on $\sum_l \oplus E^{\otimes
l}\otimes F^{\otimes (k-l)}$. In fact, let $Q_k$ be the projection
of $\mathcal{F}(n,m,u)$ onto $\sum_l \oplus E^{\otimes l}\otimes
F^{\otimes (k-l)}$, form the one-parameter unitary group $\{U_t\}$
defined by $U_t:=\sum_{k=0}^{\infty} e^{ikt}Q_k$ and set
$\gamma_t=Ad U_t$. Then $\{\gamma_t\}_{t \in \mathbb{R}}$ is a
$w^*$-continuous action of $\mathbb{R}$ on
$\mathcal{L}(\mathcal{F}(n,m,u))$ that normalizes both
$\mathcal{A}_u$ and $\mathcal{L}_u$ and
$$\Phi_k(a)=\frac{1}{2\pi}\int_0^{2\pi} e^{-ikt}\gamma_t(a)dt $$
for all $a\in \mathcal{L}(\mathcal{F}(n,m,u))$. Then $\Phi_k$
leaves $\mathcal{L}_u$ invariant.

 We can define the algebra
$\mathcal{R}_u$ generated by the right shifts $R_e$ and $R_f$
defined by $$ R_e e_{i_1}\otimes e_{i_2} \otimes \cdots \otimes
e_{i_k}\otimes f_{j_1}\otimes f_{j_2} \otimes \cdots \otimes
f_{j_l} = e_{i_1}\otimes e_{i_2} \otimes \cdots \otimes
e_{i_k}\otimes f_{j_1}\otimes f_{j_2} \otimes \cdots \otimes
f_{j_l}\otimes e $$ and $$ R_f e_{i_1}\otimes e_{i_2} \otimes
\cdots \otimes e_{i_k}\otimes f_{j_1}\otimes f_{j_2} \otimes
\cdots \otimes f_{i_l} = e_{i_1}\otimes e_{i_2} \otimes \cdots
\otimes e_{i_k}\otimes f_{j_1}\otimes f_{j_2} \otimes \cdots
\otimes f_{i_l}\otimes f.$$
 The techniques of the proof of
Proposition 2.3 of \cite{KP} can be applied here to show that the
commutant of $\mathcal{R}_u$ is $\mathcal{L}_u$. Also, mapping $
 e_{i_1}\otimes e_{i_2}
\otimes \cdots \otimes e_{i_k}\otimes f_{j_1}\otimes f_{j_2}
\otimes \cdots \otimes f_{j_l}$ to $ f_{j_l}\otimes f_{j_{l-1}}
\otimes \cdots \otimes f_{j_1}\otimes e_{i_k}\otimes e_{i_{k-1}}
\otimes \cdots \otimes e_{i_1}$, we get a unitary operator
$$W:\mathcal{F}(n,m,u) \rightarrow \mathcal{F}(n,m,u^*)$$
implementing a unitary equivalence of $\mathcal{L}_u$ and
$\mathcal{R}_{u^*}$. In fact, it is easy to check that
$R_{e_i}W=WL_{e_i}$ and $R_{f_j}W=WL_{f_j}$ for every $i,j$. To
see that the commutation relation in the range is given by $u^*$,
apply $W$ to (\ref{relations2}) to get (in the range of $W$)
$e_k\otimes f_l=\sum_{i,j} \bar{u}_{(i,j),(k,l)}f_j\otimes e_i=
\sum_{i,j} (u^*)_{(k,l),(i,j)}f_j \otimes e_i$ which is equation
(\ref{relations}) with $u^*$ instead of $u$.

As in \cite{KP}, we conclude that $(\mathcal{L}_u)'=\mathcal{R}_u$
and $(\mathcal{L}_u)''=\mathcal{L}_u$.

\end{section}

%%%%%%%%%%%%%%%%%%%%%%%%%%%%%%%%%%%%%%%%%%%%%%%%%%%%%%%%%%%%%%%%%%%%%%%%
\begin{section}{The character space and its core}

In the following proposition we describe the structure of the
character spaces $\mathcal{M}(\mathcal{L}_u)$ and
$\mathcal{M}(\mathcal{A}_u)$ (equipped with the weak$^*$
topology). Similar results were obtained in \cite{KP} for algebras
defined for higher rank graphs and in \cite{DP98} for analytic
Toeplitz algebras. (See also \cite{Po}.)

It will be convenient to write
\begin{equation}\label{variety}
V_u= \{ (z,w) \in \mathbb{C}^n \times \mathbb{C}^m \;:\; z_iw_j =
\sum_{k,l}u_{(i,j),(k,l)}z_kw_l\;\}
\end{equation}
and \begin{equation}\label{omegaU} \Omega_u=V_u\cap
(\overline{\mathbb{B}}_n \times \overline{\mathbb{B}}_m)
\end{equation} where $\mathbb{B}_n$ is the open unit ball of
$\mathbb{C}^n$.
 We refer to $V_u$ as the
\emph{variety associated with } $u$.

\begin{proposition}\label{character}
\begin{enumerate}
\item[(1)] The linear multiplicative functionals on
$\mathbb{C}[\mathbb{F}^+_u]$ are in one-to-one correspondence with
points $(z,w)$ in $V_u$.
\item[(2)] $\mathcal{M}(\mathcal{A}_u)$ is homeomorphic to
$\Omega_u$.
\item[(3)] For $(z,w)\in \Omega_u$, write $\alpha_{(z,w)}$ for the
corresponding character of $\mathcal{A}_u$. Then $\alpha_{(z,w)}$
extends to a $w^*$-continuous character on $\mathcal{L}_u$ if and
only if $(z,w)\in \mathbb{B}_n \times \mathbb{B}_m$.
\end{enumerate}
\end{proposition}
\begin{proof}
Part (1) follows immediately from (\ref{relations}). Fix $\alpha
\in \mathcal{M}(\mathcal{A}_u)$ and write $z_i=\alpha(L_{e_i})$,
$1\leq i \leq n$, and $w_i=\alpha(L_{f_j})$, $1\leq j \leq m$.
From the multiplicativity and linearity of $\alpha$ and
(\ref{relations}), it follows that $(z,w)\in V_u$. Since $\alpha$
is contractive and maps $\sum_i a_iL_{e_i}$ to $\sum_i a_i z_i$,
it follows that $\norm{z} \leq 1$ and similarly $\norm{w} \leq 1$.
Thus $(z,w)\in \Omega_u$.

For the other direction, fix first $(z,w) \in \Omega_u$ with
$\norm{z}<1$ and $\norm{w}<1$. It follows from the definition of
$\Omega_u$ and from (\ref{relations}) that $(z,w)$ defines a
linear and multiplicative map $\alpha$ on the algebra
$\mathbb{C}[\mathbb{F}^+_u]$ such that $L_{e_i}$ is mapped into
$z_i$ and $\alpha(L_{f_j})=w_j$. Abusing notation slightly, we
write $\alpha(x)$ for $\alpha(L_x)$ for every $x \in E^{\otimes
k}\otimes F^{\otimes l}$. Also, for $i=(i_1, \ldots, i_k)$ and
$j=(j_1,\ldots,j_l)$, we write $e_if_j$ for $e_{i_1}\otimes \cdots
\otimes e_{i_k}\otimes f_{j_1}\otimes \cdots \otimes f_{j_l}$.
These elements form an orthonormal basis for $E^{\otimes k}\otimes
F^{\otimes l}$ and we now set
$$ w_{\alpha}=\sum_{i,j}\sum_{k,l}
\overline{\alpha(e_if_j)}e_if_j \in \mathcal{F}(X).$$ If $p_i\geq
0$ and $p_1+\ldots +p_n=k$ then there are $\frac{k!}{p_1!\cdots
p_n!}$ terms $e_{i_1}\otimes \cdots \otimes e_{i_k}$ with
$\alpha(e_{i_1}\otimes \cdots \otimes
e_{i_k})=z_1^{p_1}z_2^{p_2}\cdots z_k^{p_k}$. It follows that
$\sum_k\sum_i |\alpha(e_i)|^2=\sum_k \sum_{i=(i_1,\ldots
,i_k)}|\alpha(e_{i_1})|^2\cdots |\alpha(e_{i_k})|^2$. Thus
 $$\norm{w_{\alpha}}^2=\sum_{i,j,k,l}|\alpha(e_if_j)|^2=
(1-\norm{z}^2)^{-1}(1-\norm{w}^2)^{-1} < \infty $$

Note that, for every $x \in E^{\otimes k}\otimes F^{\otimes l}$,
$$\langle x, w_{\alpha}\rangle = \alpha(x). $$ Thus, for $e\in E$,
$\langle x, L_e^*w_{\alpha} \rangle = \langle L_ex ,
w_{\alpha}\rangle = \alpha(e\otimes x)=\alpha(e)\alpha(x)=\langle
\alpha(e)w_{\alpha},x \rangle $ and, similarly $\langle x,L_{f}^*
w_{\alpha}\rangle = \langle \alpha(f)w_{\alpha},x\rangle$ for $f
\in F$. Thus $\langle
w_{\alpha},L_e^*w_{\alpha}\rangle=\alpha(e)\alpha(w_{\alpha})=\alpha(e)\sum
|\alpha(e_if_j)|^2=\alpha(e)\norm{w_{\alpha}}^2$. Similarly,
$\langle w_{\alpha},L_f^*w_{\alpha}\rangle=\\
\alpha(f)\alpha(w_{\alpha})=\alpha(f)\sum
|\alpha(e_if_j)|^2=\alpha(f)\norm{w_{\alpha}}^2$ for $f\in F$.
Thus if we write $\nu_{\alpha}=w_{\alpha}/\norm{w_{\alpha}}$ then
$$\alpha(x)=\langle L_x\nu_{\alpha},\nu_{\alpha}\rangle$$ for
every $x\in E^{\otimes k}\otimes F^{\otimes l}$ (for every $k,l$).
This shows that $\alpha$ is contractive and is $w^*$-continuous.
We can, therefore, extend it to an element of
$\mathcal{M}(\mathcal{L}_u)$, also denoted $\alpha$.

The analysis above shows that the image of the map $\alpha \mapsto
(z,w)\in \Omega_u$ defined above (on $\mathcal{M}(\mathcal{A}_u)$)
contains $V_u\cap (\mathbb{B}_n\times \mathbb{B}_m)$. Since
$\mathcal{M}(\mathcal{A}_u)$ is compact and the map is
$w^*$-continuous, its image contains (and, thus, is equal to)
$\Omega_u$. This completes the proof of (2). To complete the proof
of (3), we need to show that, if $(z,w)\in \Omega_u$ and the
corresponding character extends to a $w^*$-continuous character on
$\mathcal{L}_u$, then $\norm{z}<1$ and $\norm{w}<1$.

For this, write $\mathcal{L}$ for the $w^*$-closed subalgebra of
$\mathcal{L}_u$ generated by $\{L_e \;:\;e\in E\}\cup \{I\}$. Let
$P$ be the projection of $\mathcal{F}(E,F,u)$ onto
$\mathcal{F}(E)=\mathbb{C}\oplus E \oplus (E\otimes E) \oplus
\cdots $. Then $P\mathcal{L}P=P\mathcal{L}_uP$ and the map $T
\mapsto PTP$, is a $w^*$-continuous isomorphism of $\mathcal{L}$
onto $P\mathcal{L}_uP$. The latter algebra is unitarily equivalent
to the algebra $\mathcal{L}_n$ studied in \cite{DP98}. A
$w^*$-continuous character of $\mathcal{L}_u$ gives rise,
therefore, to a $w^*$-continuous character on $\mathcal{L}_n$. It
follows from \cite[Theorem 2.3]{DP98} that $z\in \mathbb{B}_n$.
Similarly, one shows that $w \in \mathbb{B}_m$.
\end{proof}
\medskip

To state the next result, we first write $u_{(i,j)}$ for the $n
\times m$ matrix whose $k,l$-entry is $u_{(i,j),(k,l)}$. Thus, the
$(i,j)$ row of $u$ provides the $n$ rows of $u_{(i,j)}$.
 We then
compute
\begin{equation}\label{computation}
\sum_{k,l} u_{(i,j),(k,l)} z_kw_l=\sum_k (\sum_l
u_{(i,j),(k,l)}w_l)z_k=\sum_k (u_{(i,j)}w)_kz_k=\langle
u_{(i,j)}w,\bar{z}\rangle. \end{equation} Write $E_{i,j}$ for the
$n\times m$ matrix  whose $i,j$-entry is $1$ and all other entries
are $0$ (so that $\langle E_{i,j}w,\bar{z}\rangle= z_iw_j$) and
write $C_{(i,j)}$ for the matrix $u_{(i,j)}-E_{i,j}$. Then the
computation above yields the following.

\begin{lemma}\label{cij} With $C_{(i,j)}$ defined as above, we
have
$$V_u=\{(z,w)\in \mathbb{C}^n\times \mathbb{C}^m\;:\;\;\;\;
\langle C_{(i,j)}w,\bar{z}\rangle =0,\;\;\;\mbox{for all}\;
i,j\}.$$
\end{lemma}

\begin{definition}\label{core}
The core of $\Omega_u$ is the subset given by
$$\Omega_u^0:=\{(z,w)\in \overline{\mathbb{B}}_n\times
\overline{\mathbb{B}}_m\;:\;\;\; C_{(i,j)}w=0,\;
C_{(i,j)}^tz=0\;\;\;\;\mbox{for all}\; i,j\}.$$
\end{definition}

 Fix $(z,w)\in \Omega_u^0$. We have $u_{(i,j)}w=E_{i,j}w$ for
all $i,j$. Thus, for every $k$, \begin{equation}\label{w}\sum_l
u_{(i,j),(k,l)}w_l=\delta_{i,k}w_j \end{equation} (where
$\delta_{i,k}$ is $1$ if $i=k$ and $0$ otherwise) and, for
$a_1,a_2, \ldots, a_n$, in $\bC$ we have  $\sum_{k,l}
u_{(i,j),(k,l)}a_kw_l=a_iw_j$. Hence, if we let $\tilde{w}^{(i)}$
be the vector in $\mathbb{C}^{mn}$ defined by
$\tilde{w}^{(i)}_{(k,l)}=\delta_{k,i}w_l$, we get
$u\tilde{w}^{(i)}=\tilde{w}^{(i)}$. Similarly, for $z$, we have
\begin{equation}\label{z}
\sum_k u_{(i,j),(k,l)}z_k=\delta_{j,l}z_i \end{equation} and for
scalars $b_1,\ldots,b_m$ we have $\sum_{k,l}
u_{(i,j),(k,l)}b_lz_k=b_jz_i$. Thus, writing $\tilde{z}_{(j)}$ for
the vector defined by $(\tilde{z}_{(j)})_{(k,l)}=\delta_{l,j}z_k$,
we have $u\tilde{z}_{(j)}=\tilde{z}_{(j)}$. The vector
$\tilde{w}^{(i)}$ in $\mathbb{C}^{nm}=\mathbb{C}^n \otimes
\mathbb{C}^m$ is also expressible as $\delta_i \otimes w$ where
$\{\delta_1,\ldots,\delta_n\}$ is the standard basis of
$\mathbb{C}^n$, and, similarly, $\tilde{z}_{(j)}=z\otimes
\delta_j$. We therefore obtain Lemma \ref{dimfix} which will be
useful in Section 6.

We note also the following companion formula. Suppose $(z,w)\in
\Omega^0_u$. Then, as we noted above,
$u\tilde{z}_{(j)}=\tilde{z}_{(j)}$ and, thus,
$u^*\tilde{z}_{(j)}=\tilde{z}_{(j)}$. Writing this explicitly, we
have, for all $i,j,l$,
\begin{equation}\label{uz}
\sum_k u_{(k,l),(i,j)}\bar{z_k}=\delta_{j,l}\bar{z_i}.
\end{equation}.

\begin{lemma}\label{dimfix}
Let $(z,w)$ be a vector in the core $\Omega_u^0$. Then
$$span\{\tilde{z}_{(j)}, \tilde{w}^{(i)}\;:\;\; 1\leq i \leq n, \;1\leq j
\leq m \}\subseteq Ker(u-I).$$ In particular,
\begin{enumerate}
\item[(i)] If the core contains a vector $(z,w)$ with $z\neq 0$,
then $dim(Ker(u-I))\geq m$.
\item[(ii)] If the core contains a vector $(z,w)$ with $w\neq 0$
then $dim(Ker(u-I))\geq n$.
\item[(iii)] If the core contains a vector $(z,w)$ with $z\neq 0$
and $w\neq 0$, then $dim(Ker(u-I))\geq m+n-1$.
\end{enumerate}
\end{lemma}

We now characterise the core in an algebraic manner in terms of
representations into the algebra $T_2$ of upper triangular
$2\times 2$ matrices. We remark that nest representations such as
these have proven useful in the algebraic structure theory of
nonself-adjoint algebra \cite{KatKri}, \cite{Sol-UcanC}.

Let $\rho : \mathbb{C}[\mathbb{F}^+_u] \to T_2$ with
$$\rho(T)=\left(
\begin{array}{cc} \rho_{1,1}(T) & \rho_{1,2}(T) \\ 0 &
\rho_{2,1}(T) \end{array} \right) $$ Then $\rho_{1,1}$ and
$\rho_{2,2}$ are characters and $\rho_{1,2}$ is a linear
functional that satisfies
\begin{equation}\label{rho}
\rho_{1,2}(TS)=\rho_{1,1}(T)\rho_{1,2}(S)+\rho_{1,2}(T)\rho_{2,2}(S)
\end{equation}
for $T,S \in \mathbb{C}[\mathbb{F}^+_u]$.

We now restrict to the case where $\rho_{1,1}=\rho_{2,2}$. By
Proposition~\ref{character}(1), both are associated with a point
$(z,w)$ in $V_u$. It follows from (\ref{rho}) that $\rho_{1,2}$ is
determined by its values on $L_{e_i}$ and $L_{f_j}$. Setting
$\lambda_i=\rho_{1,2}(L_{e_i})$ and $\mu_j=\rho_{1,2}(L_{f_j})$,
we associate with each homomorphism $\rho$ (as discussed above) a
quadruple $(z,w,\lambda,\mu)$ where $(z,w)\in V_u$ and, for every
$i,j$,
\begin{equation}
z_i\mu_j+\lambda_iw_j=\sum_{k,l} u_{(i,j),(k,l)}(w_l\lambda_k +
\mu_l z_k).
\end{equation}
(The last equation follows from (\ref{relations})). Using
(\ref{computation}) we can write the last equation as
$$\langle u_{(i,j)}w, \bar{\lambda} \rangle + \langle
u_{(i,j)}\mu,\bar{z} \rangle =z_i\mu_j+\lambda_iw_j=\langle
E_{i,j}w,\bar{\lambda}\rangle + \langle
E_{i,j}\mu,\bar{z}\rangle.$$ That is,
\begin{equation}\label{clambda}
\langle C_{(i,j)}w,\bar{\lambda}\rangle + \langle \mu,
\overline{C_{(i,j)}^tz}\rangle =0.
\end{equation}

The following lemma now follows from the definition of the core.

\begin{lemma}\label{lambdamu}
A point $(z,w)\in \Omega_u$ lies in the core $\Omega_u^0$ if and
only if every $(\lambda,\mu)\in \mathbb{C}^n\times \mathbb{C}^m$
defines a homomorphism $\rho:\mathbb{C}[\mathbb{F}^+_u]
\rightarrow T_2$ such that $$
\rho(L_{e_i})=\left(\begin{array}{cc} z_i & \lambda_i
\\ 0 & z_i \end{array}\right) $$ and
$$\rho(L_{f_j})=\left(\begin{array}{cc} w_j & \mu_j \\ 0 & w_j
\end{array}\right) $$ for all $i,j$.
\end{lemma}

\section{Automorphisms of $\cL_n$ and $ \cL_u$}

We first derive the unitary automorphisms of $\cL_n$ and $\cA_n$
associated with $U(1,n)$. These were obtained by Voiculescu
\cite{V} in the setting of the Cuntz-Toeplitz algebra. However the
automorphisms restrict to an action of $U(1,n)$ on the free
semigroup algebra. The result is rather fundamental, being a
higher dimensional version of the familiar M\"obius automorphism
group on $H^\infty$.
 For the reader's
convenience we provide complete proofs. See also the discussion in
Davidson and Pitts \cite{DP98}, and in \cite{dav-survey},
\cite{Po}.

\begin{lemma}\label{X}
Let $\alpha \in \mathbb{B}_n$ and write
\begin{enumerate}
\item[(i)] $x_0=(1-\norm{\alpha}^2)^{-1/2}$,
\item[(ii)] $\eta=x_0\alpha$, and
\item[(iii)] $X_1=(I_{\mathbb{C}^n}+\eta\eta^*)^{1/2}$.
\end{enumerate}
Then
\begin{enumerate}
\item[(1)] $\norm{\eta}^2=|x_0|^2-1$,
\item[(2)] $X_1\eta=x_0 \eta$, and
\item[(3)] $X_1^2=I+\eta\eta^*$.
\end{enumerate} In particular, the matrix $X=\left( \begin{array}{cc}
x_0 & \eta^* \\ \eta & X_1 \end{array} \right) $ satisfies
$X^*JX=J$, where $J = \left( \begin{array}{cc} 1 & 0 \\
0 & -I \end{array} \right), $
\end{lemma}
\begin{proof}
Part (1) is an easy computation and part (3) follows from the
definition of $X_1$. For (2), note that $X_1^2\eta=(I+\eta
\eta^*)\eta=\eta + \norm{\eta}^2 \eta=x_0^2 \eta$ and, for every
$\zeta \in \eta^{\perp}$, $X_1\zeta=\zeta$. Suppose
$X_1\eta=a\eta+\zeta$ ($\zeta \in \eta^{\perp}$). Then $x_0^2\eta
=X_1^2\eta=a^2\eta + \zeta$ and it follows that $a=x_0$ (as $X_1
\geq 0$) and $\zeta=0$.
\end{proof}

\vspace{3mm}

The lemma exhibits  specific matrices ($X_1$ is nonnegative) in
$U(1,n)$ associated with points in the open ball. One can
similarly check (see \cite{DP98} or \cite{Po} for example) that
the general form of a matrix
$Z$ in $U(1,n)$ is $Z=\left( \begin{array}{cc} z_0 & \eta_1^* \\
\eta_2 & Z_1
\end{array} \right) $
where
\begin{eqnarray*}
\| \eta_1 \|^2 = \| \eta_2 \|^2 = |z_0|^2 -1,\\
Z_1 \eta_1 = \bar{z_0} \eta_2 ,~~~ Z^*_1 \eta_2 = z_0 \eta_1,\\
Z^{\ast}_1 Z_1 = I_n + \eta_1 \eta^{\ast}_1, ~~~ Z_1 Z^{\ast}_1 =
I_n + \eta_2 \eta^{\ast}_2.
\end{eqnarray*}
It is these equations that are equivalent to the single matrix
equation $Z^{\ast} JZ = J$.

It is well known that the map $\theta_X$ defined on $\mathbb{B}_n$
by
$$\theta_{X}(\lambda)=\frac{X_1\lambda + {\eta}}{x_0+\langle
\lambda, {\eta}\rangle}\;\;\; ,\lambda \in \mathbb{B}_n.$$
 is an automorphism of
$\bB_n$ with inverse $\theta_{X^{-1}}$. See Lemma 4.9 of
\cite{DP98} and Lemma 8.1 of \cite{Po} for example. We make use of
this in the proof of Voiculescu's theorem below.

Let  $L_{1}, \dots , L_n$ be the generators of the norm closed
algebra $\mathcal{A}_n$ and for  $\zeta \in \mathbb{C}^n$ write
$L_{\zeta}=\sum \zeta_i L_i$. Recall that the character space
$M(\cA_n)$ is naturally identifiable with the closed ball
 $\bar{\mathbb{B}}_n$, with $\lambda$ in this ball providing a
character $\phi_\lambda$ for which $\phi_\lambda(L_i) =
\lambda_i$. The proof is a reduced version of that given above for
$M(\cA_\theta)$.

\begin{theorem}
Let $\alpha\in \mathbb{B}_n$ and let $X_1, x_0, \eta$ and $X$ be
associated with $\alpha$ as in Lemma~\ref{X}. Then

(i) there is an automorphism $\Theta_{X}$ of $\mathcal{L}_n$ such
that
\begin{equation}\label{Theta}
\Theta_{\alpha}(L_{\zeta})=(x_0I+L_{\eta})^{-1}(L_{X_1\zeta}+\langle
\zeta, \bar{\eta} \rangle I),
\end{equation}

(ii) the inverse automorphism $\Theta_{X}^{-1}$ is
$\Theta_{X^{-1}}$, and $X^{-1}$ is the matrix in $U(1,n)$
associated with $-\alpha$,

(iii) there is a unitary $U_X$ on $\cF_n$ such that for $a\in
\mathcal{A}_n$,
$$U_Xa\xi_0  =
\Theta_\alpha(a)(x_0I+L_{\eta})^{-1}\xi_0$$ and $\Theta_X(a) =
U_XaU_X^*$.

\end{theorem}

\begin{proof}
Let $\cF_n$ be the Fock space for $\cL_n$, $I_n = I_{\cF_n}$, and
let $\tilde{L} = [I_n ~~  L_1 \cdots L_n]$ viewed as an operator
from $(\bC \oplus \bC^n) \otimes \cF_n = \cF_n \oplus (\bC^n
\otimes \cF_n)$ to $\cF_n$. Then $$\tilde{L}(J \otimes I)
\tilde{L}^* = I_n - \tilde{L}\tilde{L}^* = I_n -(L_1L_1^*+\dots
L_nL_n^*) = P_0
$$
where $P_0$ is the vacuum vector projection from $\cF_n $ to
$\bC$. Also, since $XJX=J$, we have
\[
\tilde{L}(J \otimes I) \tilde{L}^* = \tilde{L}(X\otimes I_n)(J
\otimes I)(X\otimes I_n) \tilde{L}^* = [Y_0~~ Y_1](J\otimes
I)[Y_0~~ Y_1]^*
\]
where
\[
[Y_0 ~~Y_1] = [I_n ~~L]\left( \begin{array}{cc} x_0\otimes I_n &
\eta^*\otimes I_n \\
\eta\otimes I_n & X_1\otimes I_n \end{array} \right). \] Thus
$Y_0Y_0^* - Y_1Y_1^*=P_0$. Also
\[
Y_0 = x_0\otimes I_n + L(\eta \otimes I_n) = x_0I_n + L_\eta,
\]
\[
Y_1 = \eta^* \otimes I_n + L(X_1 \otimes I_n)  = \eta^*\otimes I_n
+ [L_{X_1e_1} \dots L_{X_1{e_n}}]\] where, here, $e_1, \dots ,e_n$
is the standard basis for ${\bC^n}$.

The operator $V = Y_0^{-1}Y_1$ is a row isometry $[V_1 ~\cdots ~
V_n]$, from $\bC^n \otimes \cF_n$ to $\cF_n$ with defect $1$. To
see this we compute
\[
I-VV^* = I -Y_0^{-1}Y_1Y_1^*Y_0^{*-1} = I -Y_0^{-1}(-P_0 +
Y_0Y_0^*)Y_0^{*-1}
\]
\[
= I + Y_0^{-1}P_0Y_0^{*-1} -I = \xi_0^{'}\xi_0^{'*}.
\]
Here
\[
\xi_0^{'}= Y_0^{-1}\xi_0 =
(x_0I_n+L_\eta)^{-1}\xi_0=x_0^{-1}(\sum_{j=0}^\infty(x_0^{-1}L_\eta)^j\xi_0)
\]
 and so
\[
\|\xi_0'\| = |x_0|^{-2}\sum_j|x_0|^{-2j}\|\eta\|^{2j} =
\frac{1}{x_0^2-\|\eta\|^{2}} = 1.
\]
Considering the path $t \to t\alpha$ for $0\leq t \leq 1$ and the
corresponding path of partial isometries $V$ it follows from the
stability of Fredholm index that the index of $V$ and $L$ coincide
and so in fact $V$ is a row isometry. Thus $V_1,\dots ,V_n$ are
isometries with orthogonal ranges.

We now have a contractive algebra homomorphism $\cA_n \to
\cL(\cF_n)$ determined by the correspondence $L_{e_i} \to V_i,
i=1,\dots ,n$. In fact it is an algebra endomorphism $\Theta :
\cA_n \to \cA_n$. Indeed, for $\xi = (\xi_1,\dots ,\xi_n)$ we have
\[
\Theta(L_\xi) = \sum \xi_iV_i = \sum \zeta_i Y_0^{-1}Y_1 (e_i
\otimes I_n)
\]
\[
= \sum \zeta_i(x_0I_n + L_\eta)^{-1}(\eta^*\otimes I_n +
[L_{X_1e_1} \dots L_{X_1{e_n}}])[I_n \cdots I_n]^t
\]
\[
=(x_0I_n + L_\eta)^{-1}(\langle\zeta,{\eta} \rangle  I_n +
L_{X_1\zeta}).
\]
Thus far we have followed Voiculescu's proof \cite{V}. The
following argument shows that $\Theta$ is an automorphism and is
an alternative to the calculation suggested in \cite{V}. The
calculation shows that
\[\phi_\lambda \circ \Theta_X =
\phi_{\theta_{\overline{X}}(\lambda).
}
\]
We have
\[
\phi_\lambda \circ \Theta_X(L_\zeta) = \phi_\lambda((x_0I_n +
L_\eta)^{-1}(\langle \zeta,{\eta} \rangle  I_n + L_{X_1\zeta}))
\]
\[
=(x_0 + \langle \lambda,{\eta}
\rangle)^{-1}(\langle\zeta,\eta\rangle + \langle
X_1\zeta,\overline{\lambda}\rangle) = \phi_\mu(L_\zeta)
\]
where
$$
\mu =\frac{\overline{X_1^*\overline{\lambda}} +
\overline{{\eta}}}{x_0+\langle \lambda, \overline{{\eta}}\rangle}=
\frac{\overline{X_1}{\lambda} + \overline{{\eta}}}{x_0+\langle
\lambda, \overline{{\eta}}\rangle}=\theta_{\overline{X}}(\lambda).
$$

 Write $\Theta_X$ for the contractive
endomorphism $\Theta$ of $\cA_n$ as constructed above. It follows
that the composition $\Phi = \Theta_{X^{-1}}\circ\Theta_X$ is a
contractive endomorphism which,  by the remarks preceding the
statement of the theorem,  induces the identity map on the
character space, so that $\phi_\lambda = \phi_\lambda \circ
\Phi^{-1}$ for all $\lambda \in \bB_n$. Such a map must be the
identity. Indeed,  suppose that we have the Fourier series
representation $\Phi^{-1}(L_{e_1}) = a_1L_{e_1}+\dots +a_nL_{e_n}
+ X$ where $X$ is a series with terms of total degree greater than
one. It follows that
$$
\lim_{t\to 0}~~ t^{-1}\phi_{(t,0,\dots ,0)}(\Phi^{-1}(L_{e_1})) =
a_1
$$
while
$$
\lim_{t\to 0}~~ t^{-1}\phi_{(t,0,\dots ,0)}(L_{e_1}) =
1.
$$
Since the induced map is the identity, we have $a_1 = 1 $ and $a_k
= 0$ for $k \ge 2$. In this way we see that the image of each
$L_i$ has the form $L_i+T_i$ where $T_i$ has only terms of total
degree greater than one. Since $L_i\xi_0$ is orthogonal to
$T_i\xi_0$ and $\Phi^{-1}(L_i)$ is a contraction, we have $1\geq
\Vert \Phi^{-1}(L_i)\xi_0 \Vert^2=\Vert L_i\xi_0 + T_i\xi_0
\Vert^2=\Vert L_i\xi_0 \Vert^2+ \Vert T_i\xi_0 \Vert^2=1+\Vert
T_i\xi_0 \Vert^2$. Thus $T_i\xi_0=0$ and, consequently, $T_i=0$
and so the composition $\Phi$ is the identity map.

Finally, we show that $\Theta_\alpha$ is unitarily implemented.
Define $U_X$ on $\cA_n\xi_0$ by $U_Xa\xi_0 = \Theta_X(a)\xi_0' =
\Theta_X(a)(x_0I+L_{\eta})^{-1}\xi_0$ for $a\in \cA$. Since
$\Theta_X$ is an automorphism, $(U_Xa)b\xi_0 = U_Xab\xi_0 =
\Theta_X(a)\Theta_X(b)\xi_0' = \Theta_X(a)U_Xb\xi_0,$ for $a,b \in
\cA_n$, and it follows that $U_Xa = \Theta_X(a)U_X$, as linear
transformations on the dense space $\cA_n\xi_0$.

%That $U_X$ is isometric follows from the fact proven above that
%$V=[V_1\cdots V_n]$ is a row isometry with defect 1. Indeed, this
%implies that the Wold decomposition of $V$ has a pure (nonunitary)
%part with $\xi_0'$ as a cyclic wandering vector; the
%transformation $U_X$ maps the Fock space basis element $\xi_w$
%indexed by the word $w = w(e_1,\dots ,e_n)$ to the vector
%$w(V_1,\dots ,V_n)\xi_0'$. (TO CHECKOUT :This is now a well-known
%fact, see REFERENCE?\cite{???} for example, and is also given in
%\cite{V}.) The following argument shows that in fact $U_X$ is
%unitary and that $\xi_0'$ is therefore cyclic for the Fock space.

Now, $V=[V_1,\ldots,V_n]$ is a row isometry with defect space
spanned by $\xi_0'$. The map $U_X$ maps $\xi_i=L_i\xi_0$ to
$\Theta_X(L_i)\xi_0'=V_i\xi_0'$ and, if $w=w(e_1,\ldots,e_n)$ is a
word in $e_1,\dots,e_n$ , then
$$U_X\xi_w=U_Xw(L_1,\ldots,L_n)\xi_0=
\Theta_X(w(L_1,\ldots,L_n))\xi_0'=w(V_1,\ldots,V_n)\xi_0'.
$$
 Since
$V$ is a row isometry and $\xi_0'$ is a unit wandering vector for
$V$, it follows that $\{w(V_1,\ldots,V_n)\xi_0'\}$ is an
orthonormal set. Thus, $U_X$ is an isometry. Since the range of
$U_X$ contains
$U_X\mathcal{A}_n\xi_0=\Theta_X(\mathcal{A}_n)\xi_0'=
\mathcal{A}_n(x_0I+L_{\eta})^{-1}\xi_0=\mathcal{A}_n\xi_0$ we see
that  $U_X$ is unitary.
\end{proof}\medskip

\begin{remark}
With the same calculations as in the proof above and slightly more
notation, one can show that each invertible matrix $Z \in U(1,n)$
defines an automorphism $\Theta_Z$ and that $Z\to \Theta_Z$ is an
action of $U(1,n)$ on $\cA_n$ and, in particular,
$\Theta_Z\Theta_X =\Theta_{ZX}$. Moreover, $Z \to U_Z$ is a
unitary representation of $U(1,n)$ implementing this as the
following calculation indicates.

Let $W=\left( \begin{array}{cc} w_0 & \omega^* \\
\omega & W_1
\end{array} \right) $ be the matrix in $U(1,n)$ associated with
$\beta \in \bB_n$ as in Lemma \ref{X}. Then
\[ U_XU_Wa\xi_0 = U_X(\Theta_\beta(a)(w_0+L_\omega)^{-1}\xi_0)
\]
\[= \Theta_\alpha(\Theta_\beta(a)(w_0+L_\omega)^{-1})(x_0I_n+L_\eta)^{-1}\xi_0
\]
\[
=\Theta_\alpha(\Theta_\beta(a))\Theta_\alpha((w_0+L_\omega)^{-1})(x_0I_n+L_\eta)^{-1}\xi_0
\]
\[
=\Theta_{XW}(a)\Theta_\alpha((w_0+L_\omega)^{-1})(x_0I_n+L_\eta)^{-1}\xi_0
\]
\[= \Theta_{XW}(a)\big[w_0I_n+(x_0I_n+L_\eta)^{-1}
(L_{X_1\omega}+\langle\omega,\eta\rangle I_n)
\big]^{-1}(x_0I_n+L_\eta)^{-1}\xi_0
\]
\[
=\Theta_{XW}(a)\big[w_0x_0I_n+ w_0L_\eta +
L_{X_1\omega}+\langle\omega,\eta\rangle I_n) \big]^{-1} \xi_0
\]
\[
=\Theta_{XW}(a)\big[(w_0x_0I_n + \langle\omega,\eta\rangle)I_n +
L_{\omega_0\eta+X_1\omega} \big]^{-1} \xi_0.
\]
One readily checks that this is the same as $ U_{XW}(a)\xi_0 $
\end{remark}

It is evident from the last theorem and its proof that the unitary
automorphisms of $\cA_n$ and $\cL_n$ act transitively on the open
subset $\bB_n$  associated with the weak star continuous
characters. We shall show that a version of this holds for the
unitary relation algebras with respect to the open core of the
character space. As a first step to constructing automorphisms of
$\cA_u $ we obtain unitary commutation relations for the n-tuples
$[\Theta(L_{e_1}), \dots ,\Theta(L_{e_n})]$ and $[ L_{f_1},\dots
,L_{f_m}]$ for certain automorphisms $\Theta$ of the copy of
$\cA_n$ in $\cA_u$.

\begin{lemma}\label{thetaalpha}
Suppose $(z,w)\in \Omega_u^0\cap (\mathbb{B}_n\times
\mathbb{B}_m)$. Write $\alpha$ for $\bar{z}$ and let
$\Theta:=\Theta_{\alpha}$ be as in (\ref{Theta}). Then, for every
$1\leq i \leq n$ and $1\leq j \leq m$,
\begin{equation}\label{thetau}
\Theta(L_{e_i})L_{f_j}=\sum_{k,l}u_{(i,j),(k,l)}L_{f_l}\Theta(L_{e_k}).
\end{equation}
\end{lemma}
\begin{proof}
Write $Y$ for $\eta\eta^*$ and $\beta$ for $(x_0+1)^{-1}$. Since
$X_1^2=I+\eta\eta^*$, $X_1=I+\beta \eta \eta^*=I+\beta  Y$ and
$Y=(Y_{i,j})$ where
$Y_{i,j}=\eta_i\bar{\eta}_j=x_0^2\bar{z}_iz_j$. We now compute
$$
(X_1e_i)f_j=e_if_j+\sum_t \beta Y_{t,i}e_tf_j=e_if_j+\sum_{t,k,l}
\beta  Y_{t,i} u_{(t,j),(k,l)}f_le_k
$$
$$=
\sum_{k,l}u_{(i,j),(k,l)}f_le_k +\sum_{t,k,l} \beta x_0^2
\bar{z}_tz_i u_{(t,j),(k,l)}f_le_k
$$
$$
=\sum_{k,l}u_{(i,j),(k,l)}f_le_k +\beta x_0^2 z_i\sum_{t,k,l}
\bar{z}_t u_{(t,j),(k,l)}f_le_k.
$$

 Using the core equation
(\ref{uz}), the last expression is equal to
$$\sum_{k,l}u_{(i,j),(k,l)}f_le_k +\beta x_0^2 z_i\sum_{k,l}
\delta_{j,l}\bar{z}_kf_le_k
$$
$$
=\sum_{k,l}u_{(i,j),(k,l)}f_le_k + \beta x_0^2
z_i\sum_{k}\bar{z}_kf_je_k
$$
$$=\sum_{k,l}u_{(i,j),(k,l)}f_le_k +
\beta x_0^2\sum_{k,l}(\delta_{j,l}z_i)\bar{z}_kf_le_k.
$$

 Using the
core equation (\ref{z}), this is equal to
$$\sum_{k,l}u_{(i,j),(k,l)}f_le_k + \beta x_0^2\sum_{k,l}(\sum_t
u_{(i,j),(t,l)}z_t)\bar{z}_kf_le_k
$$
$$
 =
\sum_{k,l}u_{(i,j),(k,l)}f_le_k +\\ \beta x_0^2\sum_{k,l,t}
u_{(i,j),(k,l)}z_k\bar{z}_tf_le_t
$$
$$=
\sum_{k,l}u_{(i,j),(k,l)}f_le_k + \beta \sum_{k,l,t}
u_{(i,j),(k,l)}Y_{t,k}f_le_t$$
$$= \sum_{k,l}u_{(i,j),(k,l)}f_le_k
+ \beta \sum_{k,l} u_{(i,j),(k,l)}f_lYe_k
$$
$$=
\sum_{k,l}u_{(i,j),(k,l)}f_lX_1e_k.$$
Thus
\begin{equation}\label{Xef}
L_{X_1e_i}L_{f_j}=\sum_{k,l}u_{(i,j),(k,l)}L_{f_l}L_{X_1e_k}.
\end{equation}
Next, we compute
$\sum_i\bar{z}_ie_if_j=\sum_{i,k,l}u_{(i,j),(k,l)}\bar{z}_if_le_k$.
Using (\ref{uz}), this is equal to
$\sum_{k,l}\delta_{j,l}\bar{z}_kf_le_k=\sum_k \bar{z}_kf_je_k.$
Thus
\begin{equation}
\sum_i\bar{z}_ie_if_j=\sum_i\bar{z}_if_je_i
\end{equation}
and, hence, $L_{\eta}$ commutes with $L_{f_j}$. It follows that
\begin{equation}\label{fj}
L_{f_j}(x_0I-L_{\eta})^{-1}=(x_0I-L_{\eta})^{-1}L_{f_j}.
\end{equation} We have, using (\ref{Xef}) and (\ref{fj}),
$$(x_0I-L_{\eta})^{-1}L_{X_1e_i}L_{f_j}=\sum_{k,l}u_{(i,j),(k,l)}
(x_0I-L_{\eta})^{-1}L_{f_l}L_{X_1e_k}$$
$$=
 \sum_{k,l}u_{(i,j),(k,l)}L_{f_l}(x_0I-L_{\eta})^{-1}L_{X_1e_k}.
 $$
Also, applying (\ref{z}) and (\ref{fj}), we get
$$(x_0I-L_{\eta})^{-1} \langle e_i,\eta \rangle
L_{f_j}=z_iL_{f_j}(x_0I-L_{\eta})^{-1}$$
$$=\sum_l
\delta_{j,l}z_iL_{f_l}(x_0I-L_{\eta})^{-1}$$
$$ =\sum_{k,l}u_{(i,j),(k,l)}z_kL_{f_l}(x_0I-L_{\eta})^{-1}.
$$
Subtracting the last two equations, we get (\ref{thetau}).
\end{proof}\medskip

\begin{corollary}\label{doublycom}
In the notation of Lemma~\ref{thetaalpha}, for every $i,j$,
$$L_{f_j}^*\Theta(L_{e_i})=\sum_{k,l}u_{(i,l),(k,j)}\Theta(L_{e_k})L_{f_l}^*.$$
\end{corollary}
\begin{proof}
It follows from (\ref{thetau}) that
$\Theta(L_{e_i})L_{f_l}=\sum_{k,t}
u_{(i,l),(k,t)}L_{f_t}\Theta(L_{e_k})$ for every $i,l$. Thus, for
$i,j,l$, \[L_{f_j}^*\Theta(L_{e_i})L_{f_l}L_{f_l}^*=\sum_{k,t}
u_{(i,l),(k,t)}L_{f_j}^*L_{f_t}\Theta(L_{e_k})L_{f_l}^*
\]
\[=\sum_{k,t}
u_{(i,l),(k,t)}\delta_{j,t}\Theta(L_{e_k})L_{f_l}^*= \sum_k
u_{(i,l),(k,j)}\Theta(L_{e_k})L_{f_l}^*. \]
Summing over $l$, we
get
$$ L_{f_j}^*\Theta(L_{e_i})(\sum_l L_{f_l}L_{f_l}^*)=\sum_{k,l}
u_{(i,l),(k,j)}\Theta(L_{e_k})L_{f_l}^*.$$ Now, $\sum_l
L_{f_l}L_{f_l}^*=I-P$ where $P$ is the projection onto the
subspace $\mathbb{C}\oplus E \oplus (E\otimes E) \oplus \ldots$.
Note that $P$ is left invariant under the operators in the algebra
generated by $\{L_{e_i} \;:\;1\leq i \leq n \}$ and, in
particular, by $\Theta(L_{e_i})$. Thus $
L_{f_j}^*\Theta(L_{e_i})P= L_{f_j}^*P\Theta(L_{e_i})P=0=
\sum_{k,l}u_{(i,l),(k,j)}\Theta(L_{e_k})L_{f_l}^*P$. This
completes the proof of the corollary.
\end{proof}\medskip

\begin{proposition}\label{autom}
Suppose $(z,w)\in \Omega_u^0\cap(\mathbb{B}_n\times
\mathbb{B}_m)$. Then there is a automorphism $\tilde{\Theta}_z$ of
$\mathcal{A}_u$ that is unitarily implemented and such that, for
every $X\in \mathcal{A}_u$,
\begin{equation}\label{alpha}
\alpha_{(0,w)}(\tilde{\Theta}_z^{-1}(X))=\alpha_{(z,w)}(X)
\end{equation} where $\alpha_{(z,w)}$ is the character associated
with $(z,w)$ by Proposition~\ref{character}.
\end{proposition}
\begin{proof} Let $U$ be the unitary operator implementing
$\Theta$. We can view $\mathcal{F}(n,m,u)$ as the sum
$$\mathcal{F}(n,m,u)=\sum_k F^{\otimes k} \otimes \mathcal{F}(E)
$$ where $\mathcal{F}(E)=\mathbb{C}\oplus E \oplus (E\otimes E)
\oplus \cdots$. We now let $V$ be the unitary operator whose
restriction to $F^{\otimes k} \otimes \mathcal{F}(E)$ is $I_k
\otimes U$ (where $I_k$ is the identity operator on $F^{\otimes
k}$). It is easy to check that, for every $f_j$,
$$VL_{f_j}V^*=L_{f_j}.$$ Now, fix $i$. We shall show, by
induction, that, for every $k$ and every $\xi \in F^{\otimes k}
\otimes \mathcal{F}(E)$,
\begin{equation}
(I_k \otimes U)L_{e_i}\xi = \Theta(L_{e_i})(I_k\otimes U)\xi.
\end{equation}
For $k=0$ this is just the fact that $U$ implements $\Theta$.
Suppose we know this  for $k$ and fix $f_j\in F$. Then, for $\xi
\in F^{\otimes k} \otimes \mathcal{F}(E)$ we have,
$$(I_{k+1}\otimes
U)L_{e_i}L_{f_j}\xi=\sum_{k,l}u_{(i,j),(k,l)}(I_{k+1}\otimes
U)L_{f_l}L_{e_k}\xi$$
$$ =\sum_{k,l}u_{(i,j),(k,l)}L_{f_l}(I_{k}\otimes U)L_{e_k}\xi.
$$
Applying the induction hypothesis, this is equal to
$\sum_{k,l}u_{(i,j),(k,l)}L_{f_l}\Theta(L_{e_k})(I_k\otimes
U)\xi$. Using (\ref{thetau}), this is
$\Theta(L_{e_i})L_{f_j}(I_k\otimes
U)\xi=\Theta(L_{e_i})(I_k\otimes U)L_{f_j}\xi$. Since $F^{\otimes
(k+1)}\otimes \mathcal{F}(E)$ is spanned by elements of the form
$L_{f_j}\xi$ (as above) the equality follows. From the relations
of Lemma \ref{thetaalpha} it follows that the map
$\tilde{\Theta}_z: X \to VXV^*$ defines a unitary endomorphism of
$\cA_u$. Since $\Theta$ is an automorphism of $\cA_n$ it follows
that $\tilde{\Theta}_z$ gives the desired automorphism.
\end{proof}\medskip

Clearly, in Proposition~\ref{autom}, we can interchange $z$ and
$w$ to get the following, where
${\Theta}_{z,w}=\tilde{\Theta}_z\tilde{\Theta}_w$.

\begin{proposition}\label{automw}
Suppose $(z,w)\in \Omega_u^0\cap (\mathbb{B}_n \times
\mathbb{B}_m)$. Then there is a unitary automorphism
${\Theta}_{z,w}$ of $\mathcal{L}_u$ which is a homeomorphism with
respect to the $w^*$-topologies and which restricts to an
automorphism of $\mathcal{A}_u$. Moreover, for every $X\in
\mathcal{L}_u$,
\begin{equation}\label{alpha}
\alpha_{(0,0)}({\Theta}_{z,w}^{-1}(X))=\alpha_{(z,w)}(X)
\end{equation} where $\alpha_{(z,w)}$ is the character associated
with $(z,w)$ as in  Proposition~\ref{character}.
\end{proposition}

An automorphism $\Psi$ of $\mathcal{A}_u$,  defines a map on the
character space of $\mathcal{A}_u$, namely $\phi \mapsto \phi
\circ \Psi^{-1}$. Thus using Proposition~\ref{character} we have a
homeomorphism $\theta_{\Psi}$ of $\Omega_u$. Also, since $\Omega_u
\cap (\mathbb{B}_n\times \mathbb{B}_m)$ is the interior of
$\Omega_u$, $\theta_{\Psi}$ maps $\Omega_u \cap(\mathbb{B}_n
\times \mathbb{B}_m)$ onto itself.

 Similarly, if
$\Psi$ is an automorphism of $\mathcal{L}_u$ which is a
homeomorphism with respect to the $w^*$-topologies, then
$\theta_{\Psi}$ is a homeomorphism of $\Omega_u \cap(\mathbb{B}_n
\times \mathbb{B}_m)$.

 In the following theorem we identify the relative interior of the
  core as the orbit of
 $(0,0)$ under the group of maps $\theta_{\Psi}$ associated with
 automorphisms $\Psi$.

\begin{theorem}\label{corechar}
For $(z,w)\in \mathbb{B}_n\times \mathbb{B}_m$ the following
conditions are equivalent.
\begin{enumerate}
\item[(1)] $(z,w)\in \Omega_u^0$.
\item[(2)] There exists a completely isometric automorphism
$\Psi$ of $\mathcal{L}_u$ that is a homeomorphism with respect to
the $w^*$-topologies and restricts to an automorphism of
$\mathcal{A}_u$, such that $\theta_{\Psi}(0,0)=(z,w)$.
\item[(3)] There exists an algebraic automorphism
$\Psi$ of $\mathcal{A}_u$ such that $\theta_{\Psi}(0,0)=(z,w)$.
\end{enumerate}
\end{theorem}
\begin{proof} The proof that (1) implies (2) follows from
Proposition~\ref{automw}. Clearly (2) implies (3). It is left to
show that (3) implies (1).

Given a point $(z,w)\in \Omega_u$, we saw in Lemma~\ref{lambdamu}
that, for every $(\lambda,\mu)$ satisfying (\ref{clambda}) there
is a homomorphism
$\rho_{z,w,\lambda,\mu}:\mathbb{C}[\mathbb{F}_u^+] \rightarrow
T_2$. For $(z,w) = (0,0)$ equation  (\ref{clambda}) holds for
every pair $(\lambda,\mu)$.
% For some pairs $(\lambda,\mu)$ this
%homomorphism is bounded and, thus, extends to a bounded
%homomorphism of $\mathcal{A}_u$ (into $T_2$). Let $B(\lambda,\mu)$
%be the set of all these pairs $(\lambda,\mu)$. Clearly,
%$B(0,0)=\mathbb{C}^n\times \mathbb{C}^m$ (because, in this case,
 Since $\rho_{0,0,\lambda,\mu}$ vanishes off a finite dimensional
subspace, it is a bounded homomorphism. In fact, for every
$(\lambda,\mu)$, $\norm{\rho_{0,0,\lambda,\mu}}\leq
1+\norm{\lambda}+\norm{\mu}$.

Given $\Psi$ and $(z,w)$ as in (3), for every $(\lambda,\mu)\in
\mathbb{C}^n\times \mathbb{C}^m$, $\rho_{0,0,\lambda,\mu}\circ
\Psi^{-1}$ is a homomorphism on $\mathbb{C}[\mathbb{F}_u^+]$ and,
thus, it is of the form $\rho_{z,w,\lambda',\mu'}$ for some
(unique) $(\lambda',\mu')$ satisfying (\ref{clambda}). Write
$\psi(\lambda,\mu)=(\lambda',\mu')$ and note that this defines a
continuous map. To prove the continuity, suppose
$(\lambda_n,\mu_n)\rightarrow (\lambda,\mu)$ and write $\rho_n$
for $\rho_{0,0,\lambda_n,\mu_n}$ and $\rho$ for
$\rho_{0,0,\lambda,\mu}$. Then (using the estimate on the norm of
$\rho_{0,0,\lambda,\mu}$) there is some $M$ such that
$\norm{\rho_{n}}\leq M$ for all $n$ and $\norm{\rho}\leq M$. For
every $Y\in \mathbb{C}[\mathbb{F}_u^+]$, $\rho_n(Y)\rightarrow
\rho(Y)$. Now fix $X\in \mathcal{A}_u$ and $\epsilon >0$. There is
some $Y\in \mathbb{C}[\mathbb{F}_u^+]$ such that $\norm{X-Y}\leq
\epsilon$ and there is some $N$ such that for $n\geq N$
$\norm{\rho_n(Y)-\rho(Y)}\leq \epsilon$. Thus, for such $n$,
$\norm{\rho_n(X)-\rho(X)}\leq (2M+1)\epsilon$. Setting
$X=\Psi(L_{e_i})$, we get $\lambda_n' \rightarrow \lambda'$ and
similarly for $\mu'$.

If $(z,w)$ is not in $\Omega_u^0$, then the set of all
$(\lambda,\mu)$ satisfying (\ref{clambda}) is a subspace of
$\mathbb{C}^n\times \mathbb{C}^m$ of dimension strictly smaller
than $n+m$ and, as is shown above, it contains the continuous
image (under the injective map $\psi$) of $\mathbb{C}^n\times
\mathbb{C}^m$. This is impossible.
\end{proof}\medskip

\end{section}
%%%%%%%%%%%%%%%%%%%%%%%%%%%%%%%%%%%%%%%%%%%%%%%%%%%%%

\begin{section}{Isomorphic algebras}
In this section we shall find conditions for  algebras
$\mathcal{A}_u$ and $\mathcal{A}_v$ to be (isometrically)
isomorphic. The characterisation also applies to the weak star
closed algebras $\cL_u$.

We start by considering a special type of isomorphism. We shall
now assume that the set $\{n,m\}$ for both algebras is the same.
In fact, by interchanging $E$ and $F$, we can assume that the
corresponding dimensions are the same and the algebras are defined
on $\mathcal{F}(n,m,u)$ and $\mathcal{F}(n,m,v)$ respectively.
This assumption will be in place in the discussion below up to the
end of Lemma~\ref{WAB}.

The algebra $\cA_u$ carries a natural $\mathbb{Z}_+^2$-grading,
with the $(k,l)$ labeled subspace being spanned by products of the
form $L_{e_{i_1}}L_{e_{i_2}}\dots
L_{e_{i_k}}L_{f_{i_1}}L_{f_{i_2}}\dots L_{f_{i_l}}$. Also, the
total length of such operators provides a natural
$\mathbb{Z}_+$-grading. Note that an algebra isomorphism $\Psi :
\cA_u \to \cA_v$ which respects the $\mathbb{Z}_+$-grading is
determined by a linear map between the spans of the generators\\ $
L_{e_{1}},\dots , L_{e_{n}},L_{f_{1}}, \dots ,L_{f_{m}}$. Here we
use the same notation for the generators of $\cA_u$ and $\cA_v$.
Such an isomorphism will be called \emph{graded}.

We now consider two types of graded isomorphisms, namely, either
bigraded, as in the following definition, or, in case $n=m$,
bigraded after relabeling generators.

\begin{definition}\label{bigraded}
\begin{enumerate}
\item[(i)] An isomorphism $\Psi: \mathcal{A}_u \rightarrow \mathcal{A}_v$
%(or $\Psi: \mathbb{F}^+_u \rightarrow \mathbb{F}^+_v$)
 is said to
be bigraded isomorphism if there are unitary matrices $A$
($n\times n$) and $B$ ($m\times m$) such that $$
\Psi(L_{e_i})=\sum_j a_{i,j}L_{e_j} \;,\;\;\;\Psi(L_{f_k})=\sum_l
b_{k,l}L_{f_l}.$$
\item[(ii)] If $m=n$ and $\Psi$ is a graded isomorphism such that
 $$
\Psi(L_{e_i})=\sum_j a_{i,j}L_{f_j} \;,\;\;\;\Psi(L_{f_k})=\sum_l
b_{k,l}L_{e_l}$$ for $n\times n$ unitary matrices $A$ and $B$ then
we say that $\Psi$ is a graded  exchange isomorphism.
\end{enumerate}

 We  write $\Psi_{A,B}$ for the bigraded isomorphism (as in
 (i)) and $\tilde{\Psi}_{A,B}$ for the graded  exchange
 isomorphism.
\end{definition}

Abusing notation, we write $\Psi(e_i)=\sum_j a_{i,j}e_j $ instead
of $\Psi(L_{e_i})=\sum_j a_{i,j}L_{e_j}$ for a bigraded
isomorphism (and similarly for the other expressions).

For unitary permutation matrices the following lemma was proved in
\cite[Theorem 5.1(iii)]{Po}.

\begin{lemma}\label{AB}
\begin{enumerate}
\item[(i)] If $\Psi_{A,B}$ is a bigraded isomorphism then
\begin{equation}\label{bigradedeq} (A\otimes B)v=u(A\otimes B) \end{equation}
 where
$A\otimes B$ is the $mn\times mn$ matrix whose $(i,j),(k,l)$ entry
is $a_{i,k}b_{j,l}$.
\item[(ii)] If $m=n$ and $\tilde{\Psi}_{A,B}$ is a graded
 exchange isomorphism then
\begin{equation}\label{gradedeq} (A\otimes B)\tilde{v}=u(A\otimes
B) \end{equation} where
$\tilde{v}_{(i,j),(k,l)}=\bar{v}_{(l,k),(j,i)}$.
\end{enumerate}
\end{lemma}
\begin{proof}
Assume  $\Psi=\Psi_{A,B}$ is a bigraded isomorphism. For $i,j$,
$$\Psi(e_i\otimes f_j)=(\sum_k a_{i,k}e_k)\otimes (\sum_l
b_{j,l}f_l)=\sum_{k,l}(A\otimes B)_{(i,j),(k,l)} e_k \otimes f_l
=$$ $$ \sum_{k,l,r,t}(A\otimes B)_{(i,j),(k,l)}v_{(k,l),(r,t)}f_t
\otimes e_r= \sum_{r,t} ((A\otimes B)v)_{(i,j),(r,t)}f_t\otimes
e_r.$$ On the other hand, $$\Psi(e_i\otimes
f_j)=\Psi(\sum_{k,l}u_{(i,j),(k,l)}f_l \otimes e_k)=\sum_{k,l,t,r}
u_{(i,j),(k,l)}b_{l,t}a_{k,r}f_t\otimes e_r=$$ $$ \sum_{t,r}
(u(A\otimes B))_{(i,j),(r,t)}f_t\otimes e_r.$$ This proves
equation (\ref{bigradedeq}). A similar argument can be used to
verify equation (\ref{gradedeq}).
\end{proof}\medskip

\begin{definition}\label{product}
If $u,v$ are $mn\times mn$ unitary matrices and there exist
unitary matrices $A$ and $B$ satisfying (\ref{bigradedeq}), we say
that $u$ and $v$ are product unitary equivalent.
\end{definition}

Now suppose that $A$ and $B$ are unitary matrices satisfying
(\ref{bigradedeq}). The same computation as in Lemma~\ref{AB}
shows that $W_{A,B}:E \otimes_u F \rightarrow E\otimes_v F$
defined by $$ W_{A,B}(e_i \otimes f_j)=\sum_{k,l}(A\otimes
B)_{(i,j),(k,l)} e_k \otimes f_l $$ is a well defined unitary
operator. Here the notation $E\otimes_u F$ indicates that this is
$E\otimes F$ as a subspace of $\mathcal{F}(n,m,u)$. Similarly, one
defines a unitary operator, also denoted $W_{A,B}$, from
$E^{\otimes k}\otimes F^{\otimes l}$ in $\mathcal{F}(n,m,u)$ to
$E^{\otimes k}\otimes F^{\otimes l}$ in $\mathcal{F}(n,m,v)$ by
$$W_{A,B}(e_{i_1}\otimes \cdots \otimes e_{i_k}\otimes
f_{j_1}\otimes \cdots \otimes f_{j_l})=$$ $$\sum a_{i_1,r_1}\cdots
a_{i_k,r_k}b_{j_1,t_1}\cdots b_{j_l,t_l}e_{r_1}\otimes \cdots
\otimes e_{r_k}\otimes f_{t_1}\otimes \cdots \otimes f_{t_l}.$$
This gives a well defined unitary operator
$$W_{A,B}:\mathcal{F}(n,m,u)\rightarrow \mathcal{F}(n,m,v).$$

\begin{lemma}\label{techAB}
For every $i,j$, write $Ae_i=\sum_k a_{i,k}e_k$ and $Bf_j=\sum_l
b_{j,l}f_l$. Then, for $g_1, g_2, \ldots , g_r$ in
$\{e_1,\ldots,e_n,f_1, \ldots ,f_m \}$,
\begin{equation}\label{Wg} W_{A,B}(g_1 \otimes g_2 \otimes \cdots \otimes
g_r)=Cg_1 \otimes Cg_2 \otimes \cdots \otimes Cg_r \end{equation}
where $Cg_i=Ag_i$ if $g_i \in \{e_1,\ldots,e_n\}$ and $Cg_i=Bg_i$
if $g_i \in \{f_1,\ldots,f_m\}$.
\end{lemma}
\begin{proof}
If the $g_i$'s are ordered such that the first ones are from $E$
and the following vectors are from $F$, then the result is clear
from the definition of $W_{A,B}$. Since we can get any other
arrangement by starting with one of this kind and interchanging
pairs $g_l,g_{l+1}$ successively (with $g_l\in \{e_1,\ldots,e_n\}$
and $g_{l+1}\in \{f_1,\ldots,f_m\}$), it is enough to show that
that if (\ref{Wg}) holds for a given arrangement of $e$'s and
$f$'s and we apply such an interchange, then it still holds. So,
we assume $g_l=e_k$, $g_{l+1}=f_s$ and we write $g'=g_1\otimes
\cdots \otimes g_{l-1}$, $g''=g_{l+2}\otimes \cdots \otimes g_r$,
$Cg'=Cg_1\otimes \cdots \otimes Cg_{l-1}$ and
$Cg''=Cg_{l+2}\otimes \cdots \otimes Cg_r$ and compute
%It
%is, in fact, enough to consider $g_{l+1}\otimes g_{l}$. Thus, it
%will suffice to show that (\ref{Wg}) holds for $f_l \otimes e_k$.

 $$W_{A,B}(g'\otimes f_s\otimes e_k\otimes g'')=
W_{A,B}(\sum_{i,j}\bar{u}_{(i,j),(k,s)}g'\otimes e_i \otimes
f_j\otimes g''). $$ Using our assumption, this is equal to $$
\sum_{i,j} \bar{u}_{(i,j),(k,s)} Cg'\otimes
(\sum_ta_{i,t}e_t)\otimes (\sum_q b_{j,q}f_q)\otimes Cg'' =$$ $$
\sum_{i,j,t,q}\bar{u}_{(i,j),(k,s)} a_{i,t}b_{j,q}Cg'\otimes
e_t\otimes f_q\otimes Cg''= $$
$$\sum_{i,j,t,q,d,p}\bar{u}_{(i,j),(k,s)} a_{i,t}b_{j,q}
v_{(t,q),(d,p)}Cg'\otimes f_p \otimes e_d\otimes Cg''=$$ $$\sum
(u^*)_{(k,s),(i,j)}(A\otimes
B)_{(i,j),(t,q)}v_{(t,q),(d,p)}Cg'\otimes f_p \otimes e_d\otimes
Cg'' =$$ $$ \sum_{d,p} (A\otimes B)_{(k,s),(d,p)}Cg'\otimes f_p
\otimes e_d \otimes Cg''=\sum_{d,p}a_{k,d}b_{s,p}Cg'\otimes f_p
\otimes e_d\otimes Cg'' =$$  $$Cg'\otimes Bf_s \otimes Ae_k\otimes
Cg''$$ completing the proof.
\end{proof}\medskip

The following lemma was proved in \cite[Section 7]{Po} and it
shows that the necessary conditions of Lemma~\ref{AB} are also
sufficient conditions on $A\otimes B$ for the existence of a
unitarily implemented isomorphism $\Psi_{A,B}$.

\begin{lemma}\label{WAB} For unitary matrices $A,B$ satisfying
(\ref{bigradedeq}) and $X\in \mathcal{A}_u$, the map $$ X \mapsto
W_{A,B}XW_{A,B}^*$$ is the bigraded isomorphism $\Psi_{A,B}: \cA_u
\to \cA_v$. Moreover $\Psi_{A,B}$ extends to a unitary isomorphism
$\cL_u \to \cL_v$, and similar statements holds for graded
exchange isomorphisms (when $m=n$).
\end{lemma}
\begin{proof} It will suffice to show the equality
$$\Psi_{A,B}(X)W_{A,B}=W_{A,B}X$$ for $X=L_{e_i}$ and for
$X=L_{f_j}$. Let $X=L_{f_j}$ and apply both sides of the equation
to $e_{i_1}\otimes \cdots \otimes e_{i_k}\otimes f_{j_1}\otimes
\cdots \otimes f_{j_l}$. Using Lemma~\ref{techAB}, we get
$$\Psi_{A,B}(L_{f_j})W_{A,B}(e_{i_1}\otimes \cdots \otimes
e_{i_k}\otimes f_{j_1}\otimes \cdots \otimes f_{j_l})$$ $$=\sum_r
b_{j,r}L_{f_r} ( Ae_{i_1}\otimes \cdots \otimes Ae_{i_k}\otimes
Bf_{j_1}\otimes \cdots \otimes Bf_{j_l})$$ $$=Bf_j \otimes
Ae_{i_1}\otimes \cdots \otimes Ae_{i_k}\otimes Bf_{j_1}\otimes
\cdots \otimes Bf_{j_l}$$ $$=W_{A,B}(f_j \otimes e_{i_1}\otimes
\cdots \otimes e_{i_k}\otimes f_{j_1}\otimes \cdots \otimes
f_{j_l})$$ $$=W_{A,B}L_{f_j}(e_{i_1}\otimes \cdots \otimes
e_{i_k}\otimes f_{j_1}\otimes \cdots \otimes f_{j_l}).$$ This
proves the equality for $X=L_{f_j}$. The proof for $X=L_{e_i}$ is
similar.
\end{proof}\medskip

\vspace{4mm}

At this point we drop our assumption that the set $\{n,m\}$ is the
same for both algebras and write $\{n',m'\}$ for the dimensions
associated with $\mathcal{A}_v$. We shall see in
Proposition~\ref{preserve0} (and Remark~\ref{isomLu}(i)) that, if
the algebras are isomorphic, then necessarily $\{n,m\}=\{n',m'\}$.

Given an isomorphism $\Psi:\mathcal{A}_u \rightarrow
\mathcal{A}_v$ we get a homeomorphism $\theta_{\Psi}:\Omega_u
\rightarrow \Omega_v$ (as in the discussion preceeding
Theorem~\ref{corechar}). The arguments used in the proof of
Theorem~\ref{corechar} to show that part (3)  implies part (1)
 apply also to isomorphisms and thus, $\theta_{\Psi}(0,0)\in \Omega_v^0$.

\begin{proposition}\label{thetacore}
Let $\Psi:\mathcal{A}_u \rightarrow \mathcal{A}_v$ be an
(algebraic) isomorphism. Then
$\theta_{\Psi}(\Omega_u^0)=\Omega_v^0$ and
$\theta_{\Psi}(\Omega_u^0\cap (\bB_n\times \bB_m))=\Omega_v^0\cap
(\bB_n\times \bB_m).$
\end{proposition}
\begin{proof}
 Fix $(z,w)$ in $\Omega_u^0$ and use Theorem~\ref{corechar} to
get an automorphism $\Phi$ of $\mathcal{A}_u$ such that
$\theta_{\Phi}(0,0)=(z,w)$. But then $\theta_{\Psi \circ
\Phi}(0,0)=\theta_{\Psi}(z,w)$ and, as we noted above, this
implies that $\theta_{\Psi}(z,w)\in \Omega_v^0$. It follows that
$\theta_{\Psi}(\Omega_u^0) \subseteq \Omega_v^0$ and, applying
this to $\Psi^{-1}$, the lemma follows.
\end{proof}\medskip

\begin{lemma}\label{bihol}
The map $\theta_{\Psi}$ is a biholomorphic map.
\end{lemma}
\begin{proof}
The coordinate functions for $\theta_{\Psi}$ are $(z,w)\mapsto
\alpha_{(z,w)}(\Psi^{-1}(e_i))$ (and $(z,w)\mapsto
\alpha_{(z,w)}(\Psi^{-1}(f_j))$) where $\alpha_{(z,w)}$ is the
character associated with $(z,w)$ by Proposition~\ref{character}.
For every $Y\in \mathbb{C}[\mathbb{F}_v^+]$, $\alpha_{(z,w)}(Y)$
is a polynomial in $(z,w)$ (for $(z,w)\in \Omega_v$) and,
therefore, an analytic function. Each $X\in \mathcal{A}_v$ is a
norm limit of elements in $\mathbb{C}[\mathbb{F}_v^+]$ and, thus,
$\alpha_{(z,w)}(X)$ is an analytic function being a uniform limit
of analytic functions on compact subsets of $\Omega_v$. Hence, for
every $(z,w)\in \Omega_v$, there is a power series that converges
in some, non empty, circular, neighborhood $C$ of $(z,w)$ that
represents $\alpha_{(z,w)}(X)$ on $C\cap \Omega_v$. Taking for $X$
the operators $\Psi^{-1}(e_i)$ and $\Psi^{-1}(f_j)$, we see that
$\theta$ is analytic. The same arguments apply to $\theta^{-1}$.
\end{proof}
\medskip

The facts in the following proposition obtained in \cite{Po} in
the case of permutation matrices.

\begin{proposition}\label{preserve0}
Let $\Psi: \cA_u \to \cA_v$ be an algebraic isomorphism and let
$\theta_{\Psi}: \Omega_u \to \Omega_v$ be the associated map
between the character spaces. Suppose $\theta_{\Psi}(0,0)=(0,0)$.
Then we have the following.
\begin{enumerate}
\item[(1)] $\{n,m\}=\{n',m'\}$ and we shall assume that $n=n'$ and
$m=m'$ (interchanging $E$ and $F$ and changing $u$ to $u^*$ if
necessary).
\item[(2)] There are unitary matrices
$U$ ($n\times n$) and $V$ ($m\times m$) such that
$\theta_{\Psi}(z,w)=(Uz,Vw)$ for $(z,w)\in \Omega_u$. (If $n=m$ it
is also possible that $\theta_{\Psi}(z,w)=(Vw,Uz)$.)
\item[(3)] If $\Psi$ is an isometric isomorphism,
then $\Psi$ is a bigraded isomorphism. (Or, if $m=n$, it may be
a graded exchange isomorphism).
\end{enumerate}
\end{proposition}

\begin{proof} The proof of Proposition 6.3 in \cite{Po}
giving (1) and (2) in the permutation case is based essentially on
Schwarz's lemma for holomorphic map from the unit disc. It applies
without change to the case of unitary matrices.

For (3) we may assume $m = m'$ and $n=n'$. From (2) we have for
each $\Phi(L_{e_i}) = L_{Ue_i} + X$ where $X$ is a sum of higher
order terms. Since $\Phi(L_{e_i})$ is a contraction and $L_{Ue_i}$
is an isometry it follows, as in the proof of Voiculescu's
theorem, that $X = 0$.
 Similarly, $\Phi(L_{f_j}) = L_{Vf_j}$ and it follows that $\Phi$
 is bigraded.
\end{proof}

Since every graded isomorphism $\Psi$ satisfies
$\theta_{\Psi}(0,0)=(0,0)$, we conclude the following.

\begin{corollary}\label{graded}
Every graded isometric isomorphism is bigraded if $n\neq m$ and
otherwise is either bigraded or is a graded  exchange isomorphism.
\end{corollary}

\begin{theorem}\label{bigradedisom}
The following statements are equivalent for unitary matrices $u,v$
in $M_n(\mathbb{C})\otimes M_m(\mathbb{C})$.

(i) There is an isometric isomorphism $\Psi : \mathcal{A}_u
\rightarrow \mathcal{A}_v$.

(ii) There is a graded isometric isomorphism from $\Psi :
\mathcal{A}_u \rightarrow \mathcal{A}_v$.

(iii) The matrices $u, v$ are product unitary equivalent or (in
case $n=m$) the matrices $u, \tilde{v}$ are product unitary
equivalent, where $\tilde{v}_{(i,j),(k,l)}=\bar{v}_{(l,k),(j,i)}$.

(iv) There is an isometric w*-continuous isomorphism $\Gamma :
\L_u \to \L_v$.
\end{theorem}
\begin{proof}
Given $\Psi$ in (i), let $(z,w) = \theta_\Psi(0,0)$. By
Proposition \ref{thetacore} $(z,w)$ lies in the interior of
$\Omega^0_v$. By Theorem \ref{corechar} there is a completely
isometric automorphism $\Phi$ of $\cA_v$ such that
$\theta_{\Phi}(0,0)=(z,w)$ and, therefore, $\theta_{\Phi^{-1}\circ
\Psi}(0,0)=(0,0)$. By Proposition~\ref{preserve0}, $\Phi^{-1}\circ
\Psi$
 is
a graded isometric isomorphism and (ii) holds. Lemma \ref{AB}
shows that (ii) implies (iii) and Lemma \ref{WAB} that (iii)
implies (i).

Finally, (iii) implies (iv) follows from Lemma \ref{WAB}, and (iv)
implies (ii) is entirely similar to (i) implies (ii).
\end{proof}\medskip

\begin{remark}\label{isomLu}
 The argument at the beginning of the proof of
Theorem~\ref{bigradedisom} shows that, whenever $\cA_u$ and
$\cA_v$ are isomorphic, we have $\{n,m\}=\{n',m'\}$.
\end{remark}

\begin{theorem}\label{Auautom} For $n\neq m$ the
isometric automorphisms of $\cA_u$ are of the form
$\Psi_{A,B}\Theta_{z,w}$ where $(z,w) \in \Omega_u^0$ and
$(A\otimes B)u = u(A \otimes B)$. In case $n=m$ the isometric
automorphisms include, in addition, those of the form
$\tilde{\Psi}_{A,B}\Theta_{z,w}$ where $(A\otimes B)\tilde{u} =
u(A \otimes B)$.
\end{theorem}

\end{section}

%%%%%%%%%%%%%%%%%%%%%%%%%%%%%%%%%

\begin{section}{Special cases}

\subsection{The case $n=m=2$}
%We now restrict our attention
% to the case where $n=m=2$. In this case, $u$ is a $4\times 4$
% matrix and we write $d:=dim Ker (u-I)$ and discuss separetely
% each of the $5$ possible values of $d$.
 Even in the low dimensions $n=m=2$ there are
 many isomorphism classes and special cases.
 Note that the product unitary equivalence class
orbit $\cO(u)$ of the $4\times 4$ unitary matrix $u$ takes the
form
\[
\cO(u)=\{(A\otimes B)u(A\otimes B)^*: A,B \in SU_2(\mathbb{C}) \},
\]
and so the product unitary equivalence classes are parametrised by
the set of orbits, $U_4(\mathbb{C})/Ad(SU_2(\mathbb{C})\times
SU_2(\mathbb{C}))$. This set admits a 10-fold parametrisation,
since, as is easily checked, $U_4(\mathbb{C})$ and
$SU_2(\mathbb{C})\times SU_2(\mathbb{C})$ are real algebraic
varieties of dimension $16$ and $6$ respectively. It follows that
the isometric isomorphism types of the algebras $\cA_u$ admit a
$10$ fold real parametrisation, with coincidences only for pairs
$\cO(u), \cO(v)$ with $u = \tilde{v}$

  We now look at some special cases in more detail.
  Let $d=dimKer(u-I)$.

\vspace{3mm}

 \textbf{Case I: $d=0$}

For every $(z,w)\in \overline{\mathbb{B}}_2\times
\overline{\mathbb{B}}_2$,
we have $(z,w)\in \Omega_u$ if and only if the vector\\
$(z_1w_1,z_1w_2,z_2w_1,z_2w_2)^t$ lies in $Ker(u-I)$. Thus, in
case I, $\Omega_u$ is as small as possible and is equal to $$
\Omega_{min}:=(\overline{\mathbb{B}}_2\times
\{0\})\cup(\{0\}\times \overline{\mathbb{B}}_2) .$$ It follows
from Lemma~\ref{dimfix} that, in this case,
$$\Omega_u^0=\{(0,0)\}.$$ By Proposition~\ref{preserve0} every
isometric automorphism of $\mathcal{A}_u$ is graded and the
isometric automorphisms of $\mathcal{A}_u$ are given by pairs
$(A,B)$ of unitary matrices such that $A\otimes B$ either commutes
with $u$ or intertwines $u$ and $\tilde{u}$.

\vspace{3mm}

\textbf{Case II: $d=1$}\\ When $d=1$ it still follows from
Lemma~\ref{dimfix} that $$\Omega_u^0=\{(0,0)\} $$ but now it is
possible for $\Omega_u$ to be larger than $\Omega_{min}$. In fact,
if the non zero vector $(a,b,c,d)^t$ spanning $Ker(u-I)$ satisfies
$ad\neq bc$ then $\Omega_u=\Omega_{min}$ but if $ad=bc$ then the
matrix $\left(\begin{array}{cc} a & b \\ c & d \end{array}\right)$
is of rank one and can be written as $(z_1,z_2)^t(w_1,w_2)$. Thus,
$(z,w)\in V_u$ and $\Omega_u$ contains some $(z,w)$ with non zero
$z$ and $w$.

Since $\Omega_u^0=\{(0,0)\}$, it is still true that isometric
isomorphisms and automorphisms of these algebras are graded.

\vspace{3mm}

\textbf{Case III: $d=2$}\\ When $d=2$ it is possible that
$\Omega_u^0$ will contain non zero vectors $(z,w)$ but, as
Lemma~\ref{dimfix} shows, it does not contain a vector with both
$z\neq 0$ and $w\neq 0$. All other possibilities may occur. For
example write $u_1,u_2$ and $u_3$ for the three diagonal matrices:
$$u_1=diag(1, -1, -1, 1),\;\; u_2=diag(1,-1,1,-1) $$ and
$$u_3=diag(1,1,-1,-1).$$ Using the definition of the core, we
easily see that $$\Omega_{u_1}^0=\{(0,0)\},\;\;
\Omega_{u_2}^0=\{(0,0,w_1,0):\;|w_1|\leq 1 \} $$ and $$
\Omega_{u_3}^0=\{(z_1,0,0,0):\;|z_1|\leq 1 \}. $$

Thus, the only isometric automorphisms of $\mathcal{A}_{u_1}$ are
graded, the isometric automorphisms of $\mathcal{A}_{u_2}$ are
formed by composing graded automorphisms with automorphisms of the
type described in Proposition~\ref{automw} (with $z=(0,0)$ and
$w=(w_1,0)$). Similarly, for the automorphisms of
$\mathcal{A}_{u_3}$, we use Proposition~\ref{autom}.

\vspace{3mm}

\textbf{Case IV: $d=3$}\\ In this case we are able to obtain an
explicit 2-fold parametrization of the isomorphism types of the
algebra $\mathcal{A}_u$.

 Every $4\times 4$ unitary matrix $u$
with $dim(Ker(u-I))=3$ is determined by a unit eigenvector $x$ and
its (different from $1$) eigenvalue. So that $ux=\lambda x$,
$\norm{x}=1$, $|\lambda |=1$ and $\lambda \neq 1$. Suppose $u$ and
$v$ are product unitary equivalent; that is $$(A\otimes
B)u=v(A\otimes B)$$ for unitary matrices $A,B$, and write
$x,\lambda $ for the unit eigenvector and eigenvalue of $u$. (Of
course, $x$ is determined only up to a multiple by a scalar of
absolute value $1$). Then $y=(A\otimes B)x$ is a unit eigenvector
of $v$ with eigenvalue $\lambda$. For unit vectors $x,y$ (in
$\mathbb{C}^4$) we write $x\sim y$ if there are unitary ($2\times
2$) matrices $A,B$ with $y=(A\otimes B)x$.  For the statement of
the next lemma recall that the entries of the vectors $x$ and $y$
in $\mathbb{C}^4$ are indexed by $\{(i,j): 1\leq i,j \leq 2 \}$.

\begin{lemma}\label{equiv} For a vector $x=\{x_{(i,j)}\}$ in
$\mathbb{C}^4$, write $c(x)$ for the $2\times 2$ matrix
$$c(x)=\left(\begin{array}{cc} x_{(1,1)} & x_{(1,2)} \\ x_{(2,1)}
& x_{(2,2)} \end{array}\right) .$$ Then $x\sim y$ if and only if
there are unitary matrices $A,B$ such that $c(x)=Ac(y)B$. (In this
case, we shall write $c(x)\sim c(y)$.)
\end{lemma}
\begin{proof}
Suppose $y=(A\otimes B)x$ for some unitary matrices $A=(a_{i,j}) $
and $B=(b_{i,j})$. Then $c(y)_{i,j}=y_{(i,j)}=\sum (A\otimes
B)_{(i,j),(k,l)}x_{(k,l)}=\sum_{k,l} a_{i,k}b_{j,l}c(x)_{k,l}=
(Ac(x)B)_{i,j}$.
\end{proof}\medskip

Using the polar decomposition $c(x)=U|c(x)|$ and diagonalizing
$|c(x)|=V\left( \begin{array}{cc} a & 0 \\ 0 & d \end{array}
\right) V^*$, we find that $c(x)\sim \left( \begin{array}{cc} a &
0
\\ 0 & d \end{array} \right)=c(y)$ where $y=(a,0,0,d)$ and
$a,d\geq 0$. Then $a,d$  (the eigenvalues of $|c(x)|$) are
uniquely determined once we choose them such that $a\leq d$ and,
if $\norm{x}=1$, then $a^2+d^2=1$ (so that $0\leq a \leq
1/\sqrt{2}$ and $a$ determines $d$). In this way, we associate to
each unitary matrix $u$ as above a pair $(a,\lambda)$ with $0\leq
a \leq 1/\sqrt{2}$, $\lambda \neq 1$ and $|\lambda|=1$. Using
Lemma~\ref{equiv} and the discussion preceeding it, we have the
following.

\begin{corollary}\label{produeq} For every $4\times 4$ unitary matrix $u$ with
$dim(Ker(u-I))=3$, there are numbers $\lambda$ (with $|\lambda|=1$
and $\lambda\neq 1$) and $a$ ($0\leq a \leq 1/\sqrt{2}$) such that
$u$ and $v$ are product unitary equivalent if and only if they
have the same $a,\lambda$.
\end{corollary}
\begin{proof}
Let $u$ and $v$ be unitary matrices with $dim(Ker(u-I))=3$ and let
$(a,\lambda)$, $(b,\mu)$ be the pairs associated to $u$ and $v$
(respectively) as above. Also write $x$ for the unit eigenvector
of $u$ associated to the eigenvalue $\lambda$ and let $y$ be the
unit eigenvector of $v$ associated to $\mu$.

Suppose $u$ and $v$ are product unitarily equivalent. Then they
are unitary equivalent and, thus, $\lambda=\mu$. Write $(A\otimes
B)u=v(A\otimes B)$ for unitary matrices $A,B$. As we saw above,
$y$ can be chosen to be $(A\otimes B)x$ so that $x\sim y$ and, by
Lemma~\ref{equiv}, $c(x)\sim c(y)$. It follows that $a=b$.

Conversely, assume that $a=b$ and $\lambda =\mu$. Then $c(x)\sim
c(y)$ and, thus, $x\sim y$ so we can write $y=(A\otimes B)x$ for
some unitary matrices $A,B$. Writing $v'=(A\otimes B)u(A\otimes
B)^*$, we find that $y$ is the unit eigenvector of $v'$ associated
to $\lambda$. Thus $v=v'$, completing the proof.
\end{proof}\medskip

For every $a,\lambda$ as in Corollary~\ref{produeq} we let
$u(a,\lambda)$ be the following $4\times 4$ matrix.

$$u(a,\lambda)=\left(\begin{array}{cccc} (\lambda-1)a^2+1 & 0 & 0
& (\lambda - 1)a (1-a^2)^{1/2} \\ 0 & 1 & 0 & 0 \\ 0 & 0 & 1 & 0
\\ (\lambda -1)a(1-a^2)^{1/2} & 0 & 0 & \lambda + (1-\lambda)a^2
\end{array}\right).$$ It is a straightforward computation to
verify that $dim(Ker(u-I))=3$ and that $\lambda$ is an eigenvalue
of $u(a,\lambda)$ with eigenvector $(a,0,0,(1-a^2)^{1/2})^t$. Thus
the pair associated to $u(a,\lambda)$ is $a,\lambda$ and we have

\begin{corollary}\label{ualambda}
Every matrix $u$ with $dim(Ker(u-I))=3$ is product unitary
equivalent to a unique matrix of the form $u(a,\lambda)$ (with
$0\leq a \leq 1/\sqrt{2}$, $|\lambda|=1$ and $\lambda\neq 1$).
\end{corollary}

Using the definition of the core, we immediately get the
following.

\begin{proposition}\label{core3}
If $a=0, |\lambda |=1, \lambda\ne 1$, then $\Omega_{u(0,\lambda)}$
is the union $$ \{(z_1,z_2,w_1,0):z\in \mathbb{B}_2;|w_1|\leq 1
\}\cup \{(z_1,0,w_1,w_2):w\in \mathbb{B}_2;|z_1|\leq 1 \}, $$ and
$$ \Omega_{u(0,\lambda)}^0 =\{(z_1,0,w_1,0):|z_1|\leq1;|w_1|\leq 1
\}. $$ If $a\neq 0$ then $$ \Omega_{u(a,\lambda)} =
\{(z_1,z_2,w_1,w_2):az_1w_1+(1-a^2)^{1/2}z_2w_2=0, (z,w)\in
\overline{\mathbb{B}}_2 \times \overline{\mathbb{B}}_2\}$$ and
$$\Omega_{u(a,\lambda)}^0=\{(0,0)\}.$$
\end{proposition}
\begin{proof}
The space $\Omega_{u(a,\lambda)} $ consists of points $(z,w)$ for
which
$$(z_1w_1,z_1w_2,z_2w_1,z_2w_2)^t=u(a,\lambda)(z_1w_1,z_1w_2,z_2w_1,z_2w_2)^t,$$
that is, for which
\[
((\lambda-1)a^2+1)z_1w_1+(\lambda-1)a(1-a^2)^{1/2}z_2w_2=z_1w_1,
\]
\[(\lambda-1)a(1-a^2)^{1/2}z_1w_1 + (\lambda +
(\lambda-1)a^2)z_2w_2 = z_2w_2.
\]
If $a=0$ this implies $z_2w_2 =0$, while if $a \neq 0$ then
$(z_1w_1,0,0,z_2w_2)$ is a fixed vector for $u(a,\lambda)$ and so
for some scalar $\mu$ $(z_1w_1,z_2w_2) = \mu((1-a^2)^{1/2}, -a)$.
The descriptions of $\Omega_{u(a,\lambda)}$ follows.

From the definition of the core and the fact that here $C_{12} =
C_{21}=0$ and $$C_{11} =
 \left[\begin{array}{cc} (\lambda-1)& 0\\
0&(\lambda-1)a(1-a^2)^{1/2} \end{array}\right], $$ $$C_{22} =
 \left[\begin{array}{cc} (\lambda-1)a(1-a^2)^{1/2}&0\\
 0& (\lambda-1)+(\lambda-1)a^2 \end{array}\right],
$$ we see that for $a=0$ we have $w_2=z_2=0$ while for $a\neq 0$,
$z_1=z_2=w_1=w_2=0$.
\end{proof}\medskip

Recall that, for a $4\times 4$ unitary matrix $v$ we defined the
matrix $\tilde{v}$ by
$\tilde{v}_{(i,j),(k,l)}=\bar{v}_{(l,k),(j,i)}$ and showed
(Corollary~\ref{bigradedisom}) that $\mathcal{A}_u$ and
$\mathcal{A}_v$ are isometrically isomorphic if and only if either
$u$ and $v$ or $u$ and $\tilde{v}$ are product unitary equivalent.

Now, it is easy to check that
$\widetilde{u(a,\lambda)}=u(a,\bar{\lambda})$ and so, using
Proposition~\ref{core} and previous results, we obtain the
following.

\begin{theorem}\label{isom3}
Let $0\leq a,b\leq 1/\sqrt{2}, \quad |\lambda|=|\mu|=1,\quad
\lambda,\mu \neq 1$. Then
\begin{enumerate}
\item[(1)] $\mathcal{A}_{u(a,\lambda)}$ and $\mathcal{A}_{u(b,\mu)}$ are
isometrically isomorphic if and only if $a=b$ and $\lambda$ equals
either $\mu$ or $\bar{\mu}$.
\item[(2)] When $a\neq 0$ the isometric automorphisms of
$\mathcal{A}_{u(0,\lambda)}$ are all bigraded
%(if $\lambda \neq
%-1$). If $\lambda =-1$ there is also a graded isomorphism (that is
%not bigraded).
\item[(3)] If $a=0$ then there are isometric isomorphisms
that are not graded
%(and are obtained as in Lemma~\ref{autom} or
%Lemma~\ref{automw}).
\end{enumerate}
\end{theorem}

\vspace{3mm}

\textbf{Case V: $d=4$}\\ This is the case where $u=I$. We have
$\Omega_u=\Omega_u^0=\overline{\mathbb{B}}_n \times
\overline{\mathbb{B}}_m$ and the isometric automorphisms are
obtained by composing graded automorphisms and the automorphisms
described by Proposition~\ref{autom}, Proposition~\ref{automw}.
\bigskip

\subsection{Permutation unitary relation algebras}
\medskip

With more structure assumed for a class of unitaries $u$ it may be
possible to derive an appropriately more definitive classification
of the algebras $\cA_u$. We indicate this now for the class of
permutation unitaries. A fuller discussion is in \cite{Po}.

Let $\theta \in S_4$, viewed as a permutation of the product set
$\{1,2\}\times\{1,2\} = \{11,12,21,22\}$. Associate with $\theta$
the matrix $u_\theta = u_{(i,j),(k,l)}$ where $u_{(i,j),(k,l)} =
1$ if $(k,l) = \theta(i,j)$ and is zero otherwise. If $\tau\in
S_4$ is \textit{product conjugate} to $\theta$ in the sense that
$\tau = \sigma\theta\sigma^{-1}$ with $\sigma$ in $S_2 \times
S_2$, then it follows that $u_\tau$ and $u_\theta$ are product
unitarily equivalent. Thus we need only consider product conjugacy
classes. It turns out that these classes are the same as the
product unitary equivalence classes of the matrices $u_\theta$.

It can be helpful to view a permutation $\theta$ in $S_{nm}$ as a
permutation of the entries of an $n\times m$ rectangular array,
since product conjugacy corresponds to conjugation through row
permutations and column permutations. Considering this for $n=m=2$
%and taking into account exchange isomorphisms it follows from
%Theorem \ref{bigradedisom} that there are at most 9 isomorphism
%classes of the algebras $\cA_\theta$ (and the algebras
one can verify firstly that there are at most $9$ isomorphism
types for the algebras $\cA_theta$ corresponding to the following
permutations:
\[
\theta_1 = \mbox{id}, \theta_2 = (11,12), \theta_3 = (11,22),
\]
\[ \theta_{4a}=(11,22,12), \theta_{4b} = \theta_{4a}^{-1} = (11,12,22),
\theta_5 = ((11,12),(21,22)),\]
\[\theta_6=((11,22),(12,21)), \theta_7=(11,12,22,21),
\theta_8=(11,12,21,22).
\]
The Gelfand spaces of the algebras $\cA_\theta$ (and $\cL_\theta$)
distinguish all of these algebras except for the pairs
$\{\theta_{4a}, \theta_{4b}\}$ and $\{\theta_7, \theta_8\}.$
However, one can verify in both cases that neither the pair $u,v$
nor the pair $u,\tilde{v}$ are product unitary equivalent. Theorem
\ref{bigradedisom} now applies to yield the following result from
\cite{Po}.

\begin{theorem}
For $n=m=2$ there are 9 isometric isomorphism classes for the
algebras $\cA_\theta$ and for the algebras $\cL_\theta$.
\end{theorem}

To a higher rank graph $(\Lambda, d)$ in the sense of Kumjian and
Pask \cite{KumPask} one can associate  nonself-adjoint Toeplitz
algebra $\cA_\Lambda, \cL_\Lambda$, as in Kribs and Power
\cite{KP}. In the single vertex rank $2$ case it is easy to see
that $\cA_\Lambda$ is equal  to the algebra $\cA_u$ for some
permutation matrix $u = \theta$ in $S_{nm}$. Thus Theorem
\ref{bigradedisom} classifies these algebras in terms of product
unitary equivalence restricted to $S_{nm}$ as stated formally in
the next theorem. In the rank $2$ case this is a significant
improvement on the results in \cite{Po} which, although covering
general rank, were restricted to the case of trivial core for the
character space. With $\tilde{\theta}$ the permutation for the
permutation matrix $\tilde{u}_\theta$ (which corresponds to
generator exchange) we have:

\begin{theorem}
Let $\Lambda_1$ and $\Lambda_2$ be single vertex 2-graphs with
relations determined by the permutations $\theta_1$ and
$\theta_2$. Then the rank 2 graph algebras $\cA_{\Lambda_1},
\cA_{\Lambda_2}$ are isometrically isomorphic if and only if the
pair $\theta_1,  \theta_2$ or the pair $\theta_1,
\tilde{\theta}_2$ are product unitary equivalent
\end{theorem}

 It is natural to
expect that as in the $(2,2)$ case product unitary equivalence
will correspond to product conjugacy.
\end{section}

%%%%%%%%%%%%%%%%%%%%%%%%%%%%%%%%%%%%%%%%%%%%%%%%%%%%
\begin{section}{$\mathcal{A}_u$ as a subalgebra of a tensor
algebra}

Let $\mathcal{E}_n$ be the Toeplitz extension of the Cuntz algebra
$O_n$ and write $H$ for the Fock space associated with $E$ (that
is, $H=\mathbb{C}\oplus E \oplus (E\otimes E) \oplus \cdots $).
Note that $\mathcal{E}_n$ acts naturally on $H$ ( by the ``shift"
or ``creation" operators $L_i=L_{e_i}$, $1\leq i \leq n$). In
fact, $L_{e_1}, \ldots, L_{e_n}$ generate $\mathcal{E}_n$ as a
$C^*$-algebra.

Consider also the space $\mathcal{F}(F)\otimes H=H \oplus
(F\otimes H) \oplus ((F\otimes F)\otimes H) \oplus \cdots $. This
space is isomorphic to $\mathcal{F}(E,F,u)$ and we write
$w:\mathcal{F}(F)\otimes H \rightarrow \mathcal{F}(E,F,u)$ for the
isomorphism. It will be convenient to write $w_k$ for the
restriction of $w$ to the summand $F^{\otimes k}\otimes H$ (which
is an isomorphism onto its image). Note that, for a fixed $k$,
$\{w_k^*L_{e_i}w_k \;:\;1\leq i \leq n \;\}$ is a set of $n$
isometries with orthogonal ranges. Thus it defines a
representation $\rho_k$ of $\mathcal{E}_n$ on $F^{\otimes
k}\otimes H$ (with $\rho_k(L_{e_i})=w_k^*L_{e_i}w_k$). (Note that
we are using $L_{e_i}$ for the creation operators both on $H$ and
on $\mathcal{F}(E,F,u)$. This should cause no confusion). We also
write $\rho_{\infty}$ for the representation  $\sum_k \oplus
\rho_k$ of $\mathcal{E}_n$ on $\mathcal{F}(F)\otimes H$ (where
$\rho_0$ is the representation of $\mathcal{E}_n$ on $H$).

Let $X$ be the column space $C_m(\mathcal{E}_n)$. This is a
$C^*$-module over $\mathcal{E}_n$. As a vector space it is the
direct sum of $m$ copies of $\mathcal{E}_n$. The right module
action of $\mathcal{E}_n$ on $X$ is given by $(a_i)\cdot b=(a_ib)$
and the $\mathcal{E}_n$-valued inner product is $\langle (a_i),
(b_i) \rangle =\sum_i a_i^*b_i$. For every $1\leq i \leq n$, we
write $\tilde{S}_i$ for the operator in $\mathcal{L}(X)$ defined
by
$$\tilde{S}_i (a_j)_{j=1}^m=(\sum_{j,k} u_{(i,j),(k,l)}
L_{e_k}a_j)_{l=1}^m.$$ Note that $$\langle (\sum_{j,k}
u_{(i,j),(k,l)} L_{e_k}a_j)_{l=1}^m, (\sum_{j',k'}
u_{(i,j'),(k',l)} L_{e_{k'}}b_{j'})_{l=1}^m \rangle = $$
$$ \sum_{j,j',k,k',l}\bar{u}_{(i,j),(k,l)}a_j^*L_{e_k}^*L_{e_{k'}}b_{j'}u_{(i,j'),(k',l)}
= $$ $$\sum_{j,j'} (uu^*)_{(i,j'),(i,j)}a_j^*b_{j'}=\sum_j
a_j^*b_j $$ $$ =\langle (a_j), (b_{j'}) \rangle .$$ Thus
$\tilde{S}_i$ is an isometry. A similar computation shows that
these isometries have orthogonal ranges and, thus, this family
defines a $^*$-homomorphism $\varphi:\mathcal{E}_n \rightarrow
\mathcal{L}(X)$, with $\varphi(L_{e_i})=\tilde{S}_i$, $1\leq i
\leq n$, making $X$ a $C^*$-correspondence over $\mathcal{E}_n$
(in the sense of \cite{Pi97} and \cite{MS98}). Once we have a
correspondence we can form $X \otimes X$ and, more generally,
$X^{\otimes k}$. Recall that to define $X\otimes X$ one defines
the sesquilinear form $\langle x\otimes y,x'\otimes y'
\rangle=\langle y, \varphi(\langle x,x' \rangle)y' \rangle$ on the
algebraic tensor product and then lets $X\otimes X$ be the
Hausdorff completion. The right action of $\mathcal{E}_n$ on
$X\otimes X$ is $(x\otimes y)\cdot a=x\otimes (y\cdot a)$ and the
left action is given by the map $\varphi_2$. $$
\varphi_2(a)(x\otimes y)=\varphi(a)x \otimes y.$$

 The definition of $X^{\otimes k}$ is similar (and the left action
 map is denoted $\varphi_k$)
For $k=0$ we set $X^{\otimes 0}=\mathcal{E}_n$ and $\varphi_0$ is
defined by left multiplication . Also write $\varphi_{\infty}$ for
$\sum_k \oplus \varphi_k$, the left action of $\mathcal{E}_n$ on
$\mathcal{F}(X)$.

 One can then define the Hilbert
space $X^{\otimes k}\otimes_{\mathcal{E}_n}H$ by defining the
sesquilinear form $\langle x\otimes h, y \otimes k \rangle =
\langle h, \langle x,y \rangle k \rangle$ ($x,y \in X^{\otimes
k}$) and applying the Hausdorff completion.

Now define the map $$v: X\otimes_{\mathcal{E}_n} H \rightarrow
F\otimes H $$ by setting $$v((a_i)\otimes h)=\sum_i f_i\otimes
a_ih.$$ It is straightforward to check that this map is a well
defined Hilbert space isomorphism. By induction, we also define
maps $v_{k}: X^{\otimes k}\otimes_{\mathcal{E}_n} H \rightarrow
F^{\otimes k} \otimes H$ by
\begin{equation}\label{vk}
v_{k+1}((a_j)\otimes z)=\sum_j f_j \otimes
v_k((\varphi_k(a_j)\otimes I_H)z) \end{equation} for $z\in
X^{\otimes k}\otimes_{\mathcal{E}_n} H $ and $v_0$ is the identity
map from $\mathcal{E}_n\otimes_{\mathcal{E}_n} H$ (which is
isomorphic to $H$) and $F^{\otimes 0}\otimes H =H$. Assume that
$v_k$ is a Hilbert space isomorphism of $X^{\otimes
k}\otimes_{\mathcal{E}_n} H $ onto $F^{\otimes k}\otimes H$ and
compute, for $(a_j), (b_j) \in X$ and $z,z' \in X^{\otimes
k}\otimes H$, $$\langle v_{k+1}((a_j)\otimes
z),v_{k+1}((b_j)\otimes z)\rangle= \sum_{j,j'}\langle f_j \otimes
v_k ((\varphi_k(a_j)\otimes I_H)z),f_{j'} \otimes v_k
((\varphi_k(b_{j'})\otimes I_H)z')\rangle =$$ $$ \sum_j \langle
v_k((\varphi_k(a_j)\otimes I_H)z),v_k ((\varphi_k(b_{j})\otimes
I_H)z')\rangle =$$ $$ \sum_j \langle
z,(\varphi_k(a_j^*b_{j})\otimes I_H)z')\rangle = $$ $$\langle
(a_j) \otimes z, (b_j)\otimes z' \rangle.$$  Thus, by induction,
each map $v_k$ is a Hilbert space isomorphism and, summing up, we
get a Hilbert space isomorphism
$$v_{\infty}:=\sum_k \oplus v_k
:\mathcal{F}(X)\otimes_{\mathcal{E}_n} H \rightarrow
\mathcal{F}(F)\otimes H .$$

\begin{lemma}\label{tauinf}
$v_{\infty}$ is a Hilbert space isomorphism and intertwines the
actions of $\mathcal{E}_n$. That is,
$$v_{\infty} \circ (\varphi_{\infty}(a)\otimes I_H) =
\rho_{\infty}(a) \circ v_{\infty}$$ for $a\in \mathcal{E}_n$.
\end{lemma}
\begin{proof}
We show that, for every $p\geq 0$ and $a\in \mathcal{E}_n$, we
have
\begin{equation}\label{p}
v_{p} \circ (\varphi_{p}(a)\otimes I_H) = \rho_{p}(a) \circ v_{p}.
\end{equation}
The proof will proceed by induction on $p$. For $p=0$ this is
clear so we now assume that it holds for $p$. For $1\leq i \leq
n$, $(a_j)\in X$ and $z \in X^{\otimes p}\otimes H$, we have
$v_{p+1}((\varphi_{p+1}(L_{e_i})\otimes I_H)((a_j)\otimes z))=
v_{p+1}(\varphi (L_{e_i})(a_j)\otimes z)= \sum_{l,k,j}
u_{(i,j),(k,l)} f_l \otimes v_p((\varphi_p(L_{e_k}a_j)\otimes
I_H)z)$. Using the induction hypothesis, this is equal to
$$\sum_{l,k,j}u_{(i,j),(k,l)}f_l \otimes
\rho_p(L_{e_k})\rho_p(a_j)v_pz= \sum_{l,k,j}u_{(i,j),(k,l)}f_l
\otimes w_p^* L_{e_k} w_p\rho_p(a_j)v_pz=$$
$$w_{\infty}^*\sum_{l,k,j}u_{(i,j),(k,l)}f_l \otimes e_k
\rho_p(a_j)v_pz=w_{\infty}^*\sum_j e_i \otimes f_j \otimes
\rho_p(a_j)v_pz=$$ $$\rho_{p+1}(L_{e_i})w_{p+1}^*\sum_j f_j\otimes
\rho_p(a_j)v_pz.$$ Using the induction hypothesis again, we get
$\rho_{p+1}(L_{e_i})w_{p+1}^*\sum_j f_j\otimes
v_p((\varphi_p(a_j)\otimes
I_H)z)=\rho_{p+1}(L_{e_i})v_{p+1}((a_j)\otimes z)$. This proves
(\ref{p}) for $p+1$ and the generators of $\mathcal{E}_n$. Since
both $\rho_{p+1}$ and $v_{p+1} (\varphi_{p+1}(\cdot)\otimes I_H)
v_{p+1}^*$ are $^*$-homomorphisms, (\ref{p}) holds for $p+1$ and
every $a\in \mathcal{E}_n$, completing the induction step. Thus,
(\ref{p}) holds for every $p$ and this implies the statement of
the lemma.
\end{proof}\medskip

 Write
$\delta_l$ for the vector $(a_j)$ in $X$ such that $a_l=I$ and
$a_j=0$ if $l\neq j$. The tensor algebra $\mathcal{T}_+(X)$ is
generated by the operators $T_{\delta_l}$ (where $T_{\delta_l}$ is
the creation operator on $\mathcal{F}(X)$ associated with
$\delta_l$) and the $C^*$-algebra
$\varphi_{\infty}(\mathcal{E}_n)$. The latter algebra is generated
(as a $C^*$-algebra) by the operators $\varphi_{\infty}(L_i)$
where $\{L_i\}$ is the set of generators of $\mathcal{E}_n$.

We have

\begin{lemma}\label{intertwine}
For every $1\leq i \leq n$ and $1\leq j \leq m$ and $k \geq 0$,
\begin{enumerate}
\item[(i)] $w \circ v_k \circ (\varphi_{\infty}(L_i)\otimes
I_H)=L_{e_i} \circ w \circ v_k $.
\item[(ii)] $w \circ v_{k+1} \circ (T_{\delta_j}\otimes
I_H)=L_{f_j} \circ w \circ v_k$.
\end{enumerate}
\end{lemma}
\begin{proof}
Part (i) follows from (\ref{p}) and part (ii) from (\ref{vk})
(with $\delta_j$ in place of $(a_j)$).
\end{proof}\medskip

Recalling that $w \circ v_{\infty}$ is a unitary operator mapping
$\mathcal{F}(X)\otimes H$ onto $\mathcal{F}(E,F,u)$, we get

\begin{theorem}\label{uisom}
\begin{enumerate}
\item[(1)] The algebra $\mathcal{A}_u$ is unitarily isomorphic to the
(norm closed) subalgebra of the tensor algebra $\mathcal{T}_+(X)$
that is generated by $\{\varphi_{\infty}(L_i), T_{\delta_j}\; :
\;1\leq i \leq n,\; 1\leq j \leq m \}$.
\item[(2)] The (norm closed) subalgebra of $B(\mathcal{F}(E,F,u))$
that is generated by $\{L_{e_i},L_{e_i}^*,L_{f_j}\;:\; 1\leq i
\leq n, 1\leq j \leq m\;\}$ is unitarily isomorphic to the tensor
algebra $\mathcal{T}_+(X)$ (and contains $\mathcal{A}_u$).
\item[(2)] The (norm closed) subalgebra of $B(\mathcal{F}(E,F,u))$
that is generated by $\{L_{e_i},L_{f_j}^*,L_{f_j}\;:\; 1\leq i
\leq n, 1\leq j \leq m\;\}$ is unitarily isomorphic to a tensor
algebra $\mathcal{T}_+(Y)$ (and contains $\mathcal{A}_u$).

\end{enumerate}
\end{theorem}
\begin{proof} Parts (1) and (2) follow from
Lemma~\ref{intertwine}. For part (3), note that one can
interchange the roles of $E$ and $F$. More precisely, one defines
the $C^*$-module $Y$ over $\mathcal{E}_m$ to be
$Y=C_n(\mathcal{E}_m)$ and the left action of $\mathcal{E}_m$ on
$Y$ by $\varphi_Y(L_{f_l})(b_k)_{k=1}^n=(\sum_{j,k}
\bar{u}_{(i,j),(k,l)}L_{f_j}b_k)_{i=1}^n $. This makes $Y$ into a
$C^*$-correspondence over $\mathcal{E}_m$ and the rest of the
proof proceeds along similar lines as above.
\end{proof}\medskip

Suppose $m=1$. Then $X$ is the correspondence associated with the
automorphism $\alpha$ of $\mathcal{E}_n$ given by mapping $T_i$ to
$\sum_{j=1}^n u_{i,j}T_j$ (note that $u$, in this case, is an
$n\times n$ matrix). The tensor algebra $\mathcal{T}_+(X)$ is the
analytic crossed product $\mathcal{E}_n \times_{\alpha}
\mathbb{Z}^+$ and $\mathcal{A}_u$ is unitarily isomorphic to the
subalgebra of this analytic crossed product that can be written
$\mathcal{A}_n \times_{\alpha} \mathbb{Z}^+$. One can also embed
$\mathcal{A}_u$ in $\mathcal{T}_+(Y)$ (as in
Corollary~\ref{uisom}(3)). Here $\mathcal{E}_m$ is simply the
(classical) Toeplitz algebra $\mathcal{T}$ and
$Y=C_n(\mathcal{T})$ with $\varphi_Y(T_z)(b_k)_k=(\sum_k
\bar{u}_{i,k}T_zb_k)_i$ (where $T_z$ is the generator of
$\mathcal{T}$).

\begin{remark}\label{extendauto}
Since the automorphisms $\Theta_{z,w}$ and $\Psi_{A,B}$ of $\cA_u$
are both unitarily implemented, they can be extended to
$\mathcal{T}_+(X)$. It is easy to check that they map
$\mathcal{T}_+(X)$ into itself and, thus, are automorphisms of
$\mathcal{T}_+(X)$. Hence, at least when $n\neq m$, every
automorphism of $\cA_u$ can be extended to an automorphism of the
tensor algebra $\mathcal{T}_+(X)$ that contains it (see
Theorem~\ref{Auautom}).
\end{remark}

\end{section}

\end{document}